\documentclass[aap,seceqn,citesort,MSNbibl,dvips]{arximspdf}

\doi{10.1214/09-AAP660}
\volume{20}
\issue{4}
\pubyear{2010}
\firstpage{1303}
\lastpage{1340}

\makeatletter

\newcommand{\indic}{\mathbb{I}}

\newproclaim{rem}{Remark}[section]

\newtheorem{prop}{Proposition}[section]

\newproclaim{example}{Example}[section]

\makeatother

\begin{document}
\begin{frontmatter}

\title{Lamperti-type laws}
\runtitle{Lamperti laws}

\begin{aug}
\author[A]{\fnms{Lancelot F.} \snm{James}\corref{}\thanksref{t1}\ead[label=e1]{lancelot@ust.hk}}
\runauthor{L. F. James}
\affiliation{Hong Kong University of Science and Technology}
\address[A]{Department of Information Systems\\
Hong Kong University of Science\\
\quad and Technology\\
Business Statistics and Operations Management\\
Clear Water Bay, Kowloon\\
Hong Kong\\
\printead{e1}} 
\end{aug}

\thankstext{t1}{Supported in part by Grants HIA05/06.BM03,
RGC-HKUST 6159/02P, DAG04/05.BM56 and RGC-HKUST 600907 of the HKSAR.}

\received{\smonth{8} \syear{2007}}
\revised{\smonth{11} \syear{2009}}

%
\begin{abstract}
This paper explores various distributional aspects of random
variables defined as the ratio of two independent positive random
variables where one variable has an $\alpha$-stable law, for
$0<\alpha<1$, and the other variable has the law defined by
polynomially tilting the density of an $\alpha$-stable random
variable by a factor $\theta>-\alpha$. When $\theta=0$, these
variables equate with the ratio investigated by
Lamperti [\textit{Trans. Amer. Math. Soc.}
\textbf{88} (1958) 380--387] which, remarkably, was shown to have a
simple density. This variable arises in a variety of areas and
gains importance from a close connection to the stable laws. This
rationale, and connection to the $\operatorname{PD}(\alpha,\theta)$ distribution,
motivates the investigations of its generalizations which we refer
to as Lamperti-type laws. We identify and exploit links to random
variables that commonly appear in a variety of applications.
Namely Linnik, generalized Pareto and \mbox{$z$-distributions}. In each
case we obtain new results that are of potential interest. As some
highlights, we then use these results to (i) obtain integral
representations and other identities for a class of generalized
Mittag--Leffler functions, (ii) identify explicitly the L\'evy
density of the semigroup of stable continuous state branching
processes (CSBP) and hence corresponding limiting distributions
derived in Slack and in Zolotarev [\textit{Z. Wahrsch. Verw. Gebiete}
\textbf{9} (1968) 139--145, \textit{Teor. Veroyatn. Primen.}
\textbf{2} (1957) 256--266], which are
related to the recent work by Berestycki, Berestycki and
Schweinsberg, and Bertoin and LeGall [\textit{Ann.
Inst. H. Poincar\'e Probab. Statist.} \textbf{44} (2008) 214--238,
\textit{Illinois J. Math.} \textbf{50} (2006) 147--181] on beta
coalescents. (iii) We obtain explicit results for the occupation
time of generalized Bessel bridges and some interesting stochastic
equations for $\operatorname{PD}(\alpha,\theta)$-bridges. In particular we obtain
the best known results for the density of the time spent positive
of a Bessel bridge of dimension $2-2\alpha$.
\end{abstract}

%

\end{frontmatter}

\section{Introduction}\label{sec1}

Let $S_{\alpha}$, for $0<\alpha<1$ denote a positive stable random
variable, with density $f_{\alpha}$, and having Laplace transform,
\[
\mathbb{E}[{e}^{-\lambda S_{\alpha}}]={e}^{-\lambda^{\alpha}}.
\]
Additionally, for $\theta>-\alpha$ define variables
$S_{\alpha,\theta}$ independent of $S_{\alpha}$ whose laws follow
a polynomially tilted stable distribution having density
proportional to $t^{-\theta}f_{\alpha}(t)$. When $\theta=0$,
$S_{\alpha,0}:=S'_{\alpha}\stackrel{d}{=}S_{\alpha}$. In this case
Lamperti \cite{Lamperti} (see also
Zolotarev \cite{ZolotarevStable} and Chaumont and
Yor \cite{Chaumont}) showed that, despite the general
intractability of $f_{\alpha}$, the ratio
\[
X_{\alpha}\stackrel{d}{=}\frac{S_{\alpha}}{S'_{\alpha}}
\]
has a remarkably simple density given as
%
%
\begin{equation} \label{denX}
f_{X_{\alpha}}(y)=\frac{\sin(\pi
\alpha)}{\pi}\frac{y^{\alpha-1}}{y^{2\alpha}+2y^{\alpha}\cos(\pi
\alpha)+1}\qquad{\mbox{for }}y>0.
\end{equation}

This variable arises in many important and often seemingly
unrelated contexts. For instance,
\cite{BFRY,Bourgade,DevroyeOne,DevroyeLinnik,MLP,MPS,Pillai}.
Inspired by these facts and connections to the $(\alpha,\theta)$
Poisson Dirichlet family of distributions discussed in Pitman and
Yor \cite{PY97} leads us to investigate properties of variables
defined as
\[
X_{\alpha,\theta}\stackrel{d}{=}\frac{S_{\alpha}}{S_{\alpha
,\theta}}.
\]
We refer to these variables as being \textit{Lamperti} variables or
variables having \textit{Lamperti-type} laws. Our purpose, from a
broad perspective, is to demonstrate that these variables have
strong connections to more familiar random variables that appear
in a variety of applications in probability, statistics and
related fields. In other words, the Lamperti variables, albeit
often hidden, appear in many important contexts. Furthermore, we
show how to utilize these links to both deduce properties of
$X_{\alpha,\theta}$, and develop new nontrivial results related
to the linked variables. These results can also be potentially
used to expand modeling capabilities. Our results are suggestive
of an \textit{active} beta--gamma-stable calculus that extends the
notion often associated with beta and gamma variables via
Lukacs' \cite{Lukacs} characterization.

\subsection{Outline}\label{sec11}

We now present an outline of this paper highlighting specifics.
More detailed references can be found in each section. Each
section contains new results of a nontrivial nature that in some
cases are generalizations of existing results. In addition,
combined, they represent a nice partial survey of linked
variables. Section \ref{sec2} consists of essentially two parts. The first
develops a series of pertinent distributional results for
$X_{\alpha,\theta}$ and for a broader class defined by multiplying
the Lamperti variables by beta variables. One shall notice the
class of random variables we denote as $X^{(\sigma)}_{\alpha,1}
\stackrel{d}{=}\beta_{\sigma,1-\sigma}X_{\alpha,\sigma}$ plays a
major role throughout the sections. This multiplication is based
on ideas we developed in \cite{JamesBernoulli}. The second
constitutes a natural progression of ideas, each section building
on the previous one. Specifically, Section \ref{sec23} establishes links
with positive Linnik variables. In particular, we obtain
expressions for the density of Linnik variables and also establish
an interesting gamma identity. Section \ref{sec24}, exploits this identity
in connection with generalized Pareto distributions. Albeit
brief, the main result is used to identify an unknown limiting
distribution obtained by Zolotarev \cite{ZolotarevSlack} and
Slack \cite{Slack} which we discuss in Section \ref{sec4}. Section \ref{sec24} uses the
characterization in the previous sections to demonstrate
how one can develop a calculus involving $z$-distributions. In
particular, we identify new classes of random variables, arising
as solutions to stochastic equations involving $z$-distributions,
having both hyperbolic characteristic functions and variables
whose density can be computed explicitly. Section~\ref{sec3} obtains
results for a generalization of Mittag--Leffler functions that can
be expressed as Laplace transforms of
$S^{-\alpha}_{\alpha,\theta}$ or $X_{\alpha,\theta}$ and can be
represented in terms of densities of Linnik variables. Section
\ref{sec4}
solves a fairly hard problem, identifying the explicit L\'evy
density of the semigroup of stable continuous state branching
processes. Results in Sections \ref{sec5} and \ref{sec6}, with the exception of
$\alpha=1/2$, present the best known results for occupation times
of various quantities including times spent positive on certain
random subsets. We also develop a series of interesting stochastic
equations. As one highlight we obtain results for the otherwise
elusive case of Bessel bridges of dimension $2-2\alpha$. Section
\ref{sec7}
discusses aspects of Brownian time changed models where we close
by exploiting an interesting, yet not well known, representation
of symmetric stable variables of index $0<2\alpha\leq1$, found
in \cite{DevroyePolya}.

\subsection{Some notation and background}\label{sec12}

Here we briefly recount some notation and background related to
Bessel processes and the Poisson Dirichlet family of laws. See
Pitman \cite{Pit02,Pit06} for a more precise exposition. Let
$\mathcal{B}:=(B_{t}, t>0)$ denote a strong Markov process on
$\mathbb{R}$ whose normalized ranked lengths of excursions,
$(P_{i})\in\mathcal{P}=\{\mathbf{s}=(s_{1},s_{2},\ldots):s_{1}\ge
s_{2}\ge\cdots\ge0$ and $\sum_{i=1}^{\infty}s_{i}=1\},
$ follow a Poisson Dirichlet law with parameters $(\alpha,0)$ for
$0<\alpha<1$, as discussed in Pitman and Yor \cite{PY97}. Denote
this law as $\operatorname{PD}(\alpha,0)$.
Let $(L_{t}; t>0)$ denote its local time starting at $0$, and let $\tau_{\ell}=\inf\{t\dvtx
L_{t}>\ell\}, \ell\ge0$ denote its inverse local time. In this
case $\tau$ is an $\alpha$-stable\vspace*{-2pt} subordinator where we choose
$\tau_{1}\stackrel{d}{=}S_{\alpha}$. There is the scaling identity
(see \cite{PY92}),
\[
L_{1}\stackrel{d}{=}\frac{L_{t}}{t^{\alpha}}\stackrel{d}{=}\frac
{s}{{\tau^{\alpha}_{s}}}\stackrel{d}{=}S^{-\alpha}_{\alpha},
\]
where the local time up to time $1$, $L_{1}$, satisfies
%
%
\begin{equation}\label{inverselocaltime}
L_{1}:=\Gamma(1-\alpha)^{-1}\lim_{\epsilon\rightarrow
0}\epsilon^{\alpha}|\{i\dvtx P_{i}\ge\epsilon\}| \qquad\mbox{a.s.},
\end{equation}
and is said to follow a Mittag--Leffler distribution. This shows
that $(L_{t},\tau_{t})$ have distributions determined by
$\operatorname{PD}(\alpha,0)$. Furthermore, independent of $(P_{i})$, we suppose
that for a fixed $0<p<1$, $\mathcal{B}$ is symmetrized so that
$\mathbb{P}(B_{t}>0)=p$. Under these specifications $\mathcal{B}$
could be a $p$-skewed Bessel process of dimension $2-2\alpha$. In
particular when $p=1/2,\alpha=1/2$ then $\mathcal{B}$ behaves like
a Brownian motion. An interesting aspect of $\mathcal{B}$ is the
time its spends on certain subsets of $\mathbb{R}$. Let
\[
A^{+}_{t}=\int_{0}^{t}\indic_{(B_{s}>0)}\,ds \quad\mbox{and}\quad
A^{-}_{t}=\int_{0}^{t}\indic_{(B_{s}<0)}\,ds,
\]
such that
$t=A^{+}_{t}+A^{-}_{t}$, denote the time $\mathcal{B}$ spends
positive and negative, respectively, up till time $t$. Remarkably,
by excursion theory, the time changed-processes
$(A^{+}_{\tau_{\ell}};\ell>0)$ and $(A^{-}_{\tau_{\ell}};\ell>0)$
are independent $\alpha$-stable subordinators such that
$A^{+}_{\tau_{1}}\stackrel{d}{=}p^{1/\alpha}S_{\alpha}$ and
$A^{-}_{\tau_{1}}\stackrel{d}{=}{(1-p)}^{1/\alpha}S'_{\alpha}$, for
$S'_{\alpha}\stackrel{d}{=}S_{\alpha}$. This leads to
%
%
\begin{equation}\label{subX}
X_{\alpha}\stackrel{d}{=}\frac{A^{+}_{\tau_{\ell}}}{A^{-}_{\tau
_{\ell
}}}\stackrel{d}{=}c\frac{S_{\alpha}}{S'_{\alpha}}
\quad\mbox{and}\quad
\frac{cX_{\alpha}}{cX_{\alpha}+1}\stackrel{d}{=}\frac{A^{+}_{\tau
_{\ell
}}}{\tau_{\ell}}\stackrel{d}{=}A^{+}_{1}.
\end{equation}
Hereafter, denote the law that governs $\mathcal{B}$ and related
functionals under the above specifications as
$\mathbb{P}^{(p)}_{\alpha,0}$. Denote the\vspace*{-1pt} corresponding
expectation operator as $\mathbb{E}^{(p)}_{\alpha,0}$. Thus
writing $\mathbb{P}^{(p)}_{\alpha,0}(A^{+}_{1}\in dx)/dx$ equates
with the density of the time spent positive on $[0,1]$ of a
$p$-skewed Bessel process. As noticed by Barlow, Pitman and Yor \cite
{BPY} and Pitman and Yor \cite{PY92}, this
law was originally obtained by Lamperti \cite{Lamperti}, and from
({\ref{subX}}) it equates with
\[
\mathbb{P}^{(p)}_{\alpha,0}(A^{+}_{1}\in
dx)/dx=\mathbb{P}\bigl(cX_{\alpha}/(cX_{\alpha}+1)\in dx\bigr)/dx.
\]
Now for $\theta>-\alpha$ let $\mathbb{P}^{(p)}_{\alpha,\theta}$ and
$\mathbb{E}^{(p)}_{\alpha,\theta}$ denote the law and expectation
operator of functionals connected to a $p$-skewed process whose
excursion lengths follow $\operatorname{PD}(\alpha,\theta)$. In particular if
$(P_{i})$ is distributed according to
$\operatorname{PD}(\alpha,\theta)$, then it satisfies, for measurable
$H$,
\[
\mathbb{E}^{(p)}_{\alpha,\theta}[H((P_{i}))]=\frac{\Gamma(\theta
+1)}{\Gamma(\theta/\alpha+1)}\mathbb{E}^{(p)}_{\alpha
,0}[H((P_{i}))\tau
^{-\theta}_{1}].
\]
When $\theta=\alpha$, this corresponds to the case of a Bessel
bridge of dimension $2-2\alpha$. We use the notation
$A^{(br)}_{1}$ for the variable that satisfies
\[
\mathbb{P}^{(p)}_{\alpha,\theta}\bigl(A^{(br)}_{1}\in
dx\bigr)=\mathbb{P}^{(p)}_{\alpha,\theta+\alpha}(A^{+}_{1}\in dx).
\]
Note also that under $\mathbb{P}^{(p)}_{\alpha,\theta}$,
$L_{1}\stackrel{d}{=}S^{-\alpha}_{\alpha,\theta}$, which is also
equivalent in distribution to the $\alpha$-diversity of a
$\operatorname{PD}(\alpha,\theta)$ law.

Let now $U_{1},U_{2},\ldots,$ denote a sequence of i.i.d.
uniform $[0,1]$ random variables for $(P_{i})$ distributed
according to $\operatorname{PD}(\alpha,\theta)$, and $0\leq u\leq1$, the class of
$\operatorname{PD}(\alpha,\theta)$ random cumulative distribution functions are
defined as
\[
P_{\alpha,\theta}(u)\stackrel{d}{=}\sum_{k=1}^{\infty}P_{k}\indic
_{(U_{k}\leq
u)}.
\]
Furthermore under $\mathbb{P}^{(u)}_{\alpha,\theta}$,
\[
A^{+}_{1}\stackrel{d}{=}P_{\alpha,\theta}(u)
\]
for a fixed $u$. See Bertoin \cite{BerFrag} for applications to
coagulation/fragmentation phenomena where it is called a
$\operatorname{PD}(\alpha,\theta)$-bridge and Ishwaran and
James \cite{IJ2001} (see also Pitman \cite{Pit96}) for applications
to Bayesian statistics where in particular $P_{\alpha,\theta}$ is
referred to as a Pitman--Yor process. Under this name the process
has also been applied to problems arising in natural language
processing (see Teh \cite{Teh}). When $\theta>0$ and $\alpha=0$
$P_{0,\theta}$ is a Dirichlet process which has, since the seminal
work of Ferguson \cite{Ferg73}, played a fundamental role in
Bayesian nonparametric statistics and related areas.
\begin{rem}
The $P_{\alpha,\theta}$ processes can be defined more generally by
replacing $(U_{k})$ with i.i.d. random variables $(R_{k})$ having
common distribution $F_{R}$.
\end{rem}
\begin{rem}
For basic notation, we write $\gamma_{a}$ and $\beta_{a,b}$ to
denote a gamma random variable with shape $a$ and scale $1$, and a
beta random variable with parameters $(a,b)$. If $a>0$, and $b=0$,
then we use $\beta_{a,0}:=\lim_{b\rightarrow0}\beta_{a,b}=1$.
Additionally $\xi_{\sigma}$ denotes a Bernoulli variable with
success parameter $0<\sigma\leq1$. If $X$ and $Y$ are random
variables we will assume that $XY$ is a product of independent
random variables unless otherwise specified, if we write $X, X'$
this will mean $X\stackrel{d}{=}X'$ but they are not equal. Last we
always consider $c^{\alpha}=p/q=p/(1-p)$ where $q=1-p$ unless
otherwise specified.
\end{rem}

\section{Distributional results for $X_{\alpha,\theta}$}\label{sec2}

In this section we shall derive various distributional properties of
$X_{\alpha,\theta}$. For $\tau>0$ and $0<\sigma\leq1$, we will
sometimes work with the parametrization $\tau\sigma$, to
accommodate values such as $\tau\sigma=\theta>0$ and $\tau\sigma
=\theta+\alpha$. First, we briefly discuss some pertinent
properties of random variables referred to as Dirichlet means and
the related class of infinitely divisible random variables whose
distributions are generalized gamma convolutions (GGC), as they
will play a significant role in our exposition. For more details
and related notions, one may consult \cite
{BondBook,CifarelliRegazzini79,CifarelliRegazzini,CifarelliMelilli,Ethier,JLP,JRY,LijoiMean}
and, in particular for this exposition,~\cite{JamesBernoulli}.

For a generic positive random variable $M$, let
\[
\mathcal{C}_{\tau\sigma}(\lambda;M)=\mathbb{E}[(1+\lambda
M)^{-\tau\sigma}]=\mathbb{E}[e^{-\lambda
\gamma_{\tau\sigma}M}]
\]
denote its Cauchy--Stieltjes transform of
order $\tau\sigma$. Similar to Laplace transforms,
$\mathcal{C}_{\tau\sigma}(\lambda;M)$ uniquely characterizes the
law of $M$. Let $R$ denote a nonnegative random variable with
distribution function $F_{R}$. A random variable $M$, depending on
parameters $(\tau\sigma,{R})$, is said to be a Dirichlet mean of
order $\tau\sigma$ if
%
%
\begin{equation}\label{loglaplace}
-\log\mathcal{C}_{\tau\sigma}(\lambda; M)=\tau\sigma
\mathbb{E}[\log(1+\lambda R)]:=\tau\sigma\psi_{R}(\lambda)<\infty.
\end{equation}
Equivalently $M$ satisfies the stochastic equation
\[
M\stackrel{d}{=}\beta_{\tau\sigma,1}M+(1-\beta_{\tau\sigma,1})R.
\]
We
denote such variables as $M\stackrel{d}{=}M_{\tau\sigma}(F_{R})$.

Importantly, Cifarelli and
Regazzini \cite{CifarelliRegazzini} (see
also \cite{CifarelliMelilli}), apply an inversion formula to obtain
an expression for the distribution of $M_{\tau\sigma}(F_{R})$. In
general these are expressed in terms of Abel-type transforms. An
exception is the case of $\tau\sigma=1$, where the density of
$M_{1}(F_{R})$ can be expressed as
%
%
\begin{equation} \label{M1}
\frac{1}{\pi}\sin(\pi F_{R}(x))e^{-\Phi_{R}(x)},
\end{equation}
where
$ \Phi_{R}(x)=\mathbb{E}[{\log}|x-R|\indic_{(R\neq x)}]$.
Additionally,
\[
\tau\sigma\psi_{R}(\lambda)=\tau\sigma\int_{0}^{\infty
}(1-e^{s\lambda})s^{-1}\mathbb{E}[e^{-s/R}]\,ds
\]
is also the L\'evy exponent of an infinitely divisible random
variable with L\'evy density $\tau\sigma s^{-1}\mathbb{E}[e
^{-s/R}]$. We say that such a random variable is
$\operatorname{GGC}(\tau\sigma,R)$ and may be represented in distribution as a
gamma scale mixture
%
%
\begin{equation}
\label{GGCmix}
\gamma_{\tau\sigma}M_{\tau\sigma}(F_{R})=\gamma_{\tau}\beta
_{\tau\sigma
,\tau(1-\sigma)}M_{\tau\sigma}(F_{R}).
\end{equation}
Highly relevant to (\ref{GGCmix}), and our exposition, is a result
by James \cite{JamesBernoulli}, that for each $0<\sigma\leq1$,
%
%
\begin{equation}
\label{betascaling}
\beta_{\tau\sigma,\tau(1-\sigma)}M_{\tau\sigma}(F_{R})=M_{\tau
}(F_{R\xi
_{\sigma}}),
\end{equation}
where $\xi_{\sigma}$ is\vspace*{-2pt} a Bernoulli variable with success
probability $\sigma$. Note also that
$\beta_{\tau\sigma,\tau(1-\sigma)}\stackrel{d}{=}M_{\tau}(F_{\xi
_{\sigma}})$.
One consequence is that a $\operatorname{GGC}(\tau\sigma,R)$ variable is also a
$\operatorname{GGC}(\tau,R\xi_{\sigma})$ variable. In other words, for a fixed
$\theta>0$, a $\operatorname{GGC}(\theta, R)$ variable is a
$\operatorname{GGC}(\theta',R\xi_{\theta/\theta'})$ variable for all
$\theta'>\theta$. As pointed out in \cite{JamesBernoulli}, one
significant point about these multiple representations is that if
$0<\theta=\sigma\leq1$, then one can set $\theta'=1$ and use the
explicit density formula for Dirichlet means of order $1$,
(\ref{M1}), established by Cifarelli and
Regazzini \cite{CifarelliRegazzini} to obtain an explicit
representation of the density of such a $\operatorname{GGC}(\sigma,R)$ variable.
See \cite{JamesBernoulli} for its precise form and further
details.
\begin{rem}
Letting $F^{-1}_{R}$ denote a quantile function, variables
$M_{\tau\sigma}(F_{R})$ are called Dirichlet means since they can
always be represented as
\[
M_{\tau\sigma}(F_{R})\stackrel{d}{=}\int
_{0}^{1}F^{-1}_{R}(u)P_{0,\tau
\sigma}(du)\stackrel{d}{=}\int_{0}^{\infty}yD_{\tau\sigma}(dy),
\]
where
$D_{\tau\sigma}(y)\stackrel{d}{=}\sum_{k=1}^{\infty}P_{k}\indic
_{(R_{k}\leq
y)}$ is a Dirichlet process with $(P_{k})\sim
\operatorname{PD}(0,\tau\sigma)$ and where $(R_{k})$ are i.i.d. $F_{R}$.
\end{rem}

\subsection{Identities}\label{sec21}

For the case of $X_{\alpha,\theta}$, one can show that for
$\theta>0$,
%
%
\begin{equation}
\label{Cauchy1}
\mathcal{C}_{\theta}(\lambda;X_{\alpha,\theta})=(1+\lambda
^{\alpha
})^{-{\theta}/{\alpha}}=\mathbb{E}[e^{-\lambda\gamma_{\theta}X_{\alpha,\theta}}],
\end{equation}
and
for $\theta>-\alpha$,
%
%
\begin{equation}
\label{Cauchy2}\hspace*{28pt}
\mathcal{C}_{1+\theta}(\lambda;X_{\alpha,\theta})=\mathbb
{E}[e^{-\lambda
\gamma_{\theta+1}X_{\alpha,\theta}}]=(1+\lambda^{\alpha})^{-
({\theta
+\alpha})/{\alpha}}=
\mathcal{C}_{\theta+\alpha}(\lambda;X_{\alpha,\theta+\alpha}).
\end{equation}
We will use (\ref{Cauchy1}) and (\ref{Cauchy2}) to more easily
establish the next series of results. However, we note that the
expressions in (\ref{Cauchy1}) and (\ref{Cauchy2}) are not obvious. We
will provide justification for (\ref{Cauchy1}) when we discuss
Linnik variables in the next section. Assuming that
(\ref{Cauchy1}) is true, (\ref{Cauchy2}) then follows from an
identity due to Perman, Pitman and Yor \cite{PPY92},
%
%
\begin{equation}\label{keyid1}
\frac{1}{S_{\alpha,\theta}}\stackrel{d}{=}\frac{\beta_{\theta
+\alpha
,1-\alpha}}{S_{\alpha,\theta+\alpha}}
\end{equation}
for $\theta>-\alpha$. (\ref{keyid1}) is another highly relevant
component to our exposition and shows that
$X_{\alpha,\theta}\stackrel{d}{=}\beta_{\theta+\alpha,1-\alpha
}X_{\alpha
,\theta+\alpha}$.
\begin{prop}\label{prop21}
The random variables $X_{\alpha,\theta}$ are Dirichlet
means having the
following properties: for $\theta>0$,
%
%
\begin{equation}\label{Xid1}
X_{\alpha,\theta}\stackrel{d}{=}\beta_{\theta,1}X_{\alpha,\theta
}+(1-\beta_{\theta,1})X_{\alpha}\stackrel{d}{=}M_{\theta
}(F_{X_{\alpha}}),
\end{equation}
and for $\theta>-\alpha$ and $\sigma=(\theta+\alpha)/(1+\theta)$,
%
%
\begin{equation} \label{Xid2}\qquad
X_{\alpha,\theta}\stackrel{d}{=}\beta_{\theta+\alpha,1-\alpha
}X_{\alpha
,\theta+\alpha}\stackrel{d}{=}\beta^{1/\alpha}_{(({\theta
+\alpha
})/{\alpha},({1-\alpha})/{\alpha})}X_{\alpha,1+\theta}
=M_{1+\theta}(F_{X_{\alpha}\xi_{\sigma}})
\end{equation}
with $
X_{\alpha,\theta}\stackrel{d}{=}\beta_{1+\theta,1}X_{\alpha
,\theta
}+(1-\beta_{1+\theta,1})X_{\alpha}\xi_{\sigma}$.
As special cases of (\ref{Xid2}):
\begin{longlist}
\item $X_{\alpha,1}\stackrel{d}{=}\beta_{1+\alpha,1-\alpha
}X_{\alpha
,1+\alpha};$
\item
$
X_{\alpha,1-\alpha}\stackrel{d}{=}\beta_{1,1-\alpha}X_{\alpha
,1}\stackrel{d}{=}\beta^{1/\alpha}_{({1}/{\alpha},
({1-\alpha
})/{\alpha})}X_{\alpha,2-\alpha};
$
\item
$X_{\alpha}\stackrel{d}{=}\beta_{\alpha,1-\alpha}X_{\alpha,\alpha
}\stackrel{d}{=}\beta^{1/\alpha}_{(1,({1-\alpha})/{\alpha
})}X_{\alpha,1}
=M_{1}(F_{X_{\alpha}\xi_{\alpha}})
$
which yields the identity,
\[
X_{\alpha}\stackrel{d}{=}UX_{\alpha}+(1-U)X'_{\alpha}\xi_{\alpha}
\]
for $X'_{\alpha}\stackrel{d}{=}X_{\alpha}$.
\end{longlist}
\end{prop}
\begin{pf}
In order to establish (\ref{Xid1}), we will calculate
the Cauchy--Stieltjes transform of order $\theta+1$ of the
variables appearing on the two sides of the first equality. This
entails multiplication by an independent $\gamma_{\theta+1}$
variable. Hence (\ref{Xid1}) is true if
\[
\gamma_{\theta+1}X_{\alpha,\theta}\stackrel{d}{=}\gamma_{\theta
}X_{\alpha,\theta}+\gamma_{1}X_{\alpha}.
\]
Applications of (\ref{Cauchy1}) and (\ref{Cauchy2}) show that,
$\mathcal{C}_{1+\theta}(\lambda;X_{\alpha,\theta})=\mathcal
{C}_{\theta
}(\lambda;X_{\alpha,\theta})\mathcal{C}_{1}(\lambda;X_{\alpha})$,
concluding the result. We will use similar arguments elsewhere but
will omit such details. For (\ref{Xid2}), we again calculate
$\mathcal{C}_{1+\theta}(\lambda;X_{\alpha,\theta})$. The first
equality is easily checked. For the second we use
%
%
\begin{equation}
\label{powerid}\hspace*{28pt}
\mathcal{C}_{1+\theta}\bigl(\lambda;\beta^{1/\alpha}_{(({\theta
+\alpha
})/{\alpha},({1-\alpha})/{\alpha})}X_{\alpha,1+\theta}\bigr)=\mathcal{C}_{
({1+\theta})/{\alpha}}\bigl(\lambda^{\alpha},\beta_{(({\theta
+\alpha
})/{\alpha},({1-\alpha})/{\alpha})}\bigr).
\end{equation}
In order\vspace*{-2pt} to establish the equivalence to
$M_{1+\theta}(F_{X_{\alpha}\xi_{\sigma}})$, first note that
(\ref{Xid1}) establishes
$X_{\alpha,\theta+\alpha}\stackrel{d}{=}M_{\theta+\alpha
}(F_{X_{\alpha}})$.
The result is then concluded by an application of
(\ref{betascaling}) for $\tau=1+\theta$,
$\sigma=(\theta+\alpha)/(1+\theta)$, and $R=X_{\alpha}$.
\end{pf}

The next result establishes results for the larger class of
variables defined with (\ref{GGCmix}) and (\ref{betascaling}) in
mind, as
\[
X^{(\sigma)}_{\alpha,\tau}\stackrel{d}{=}\beta_{\tau\sigma,\tau
(1-\sigma
)}X_{\alpha,\tau\sigma}.
\]
Equation (\ref{Xid2}) of Proposition \ref{prop21} is an important special
case.
\begin{prop}\label{prop22}
For $\tau>0$ and $0<\sigma\leq1$, the random variables
$X^{(\sigma)}_{\alpha,\tau}\stackrel{d}{=}\beta_{\tau\sigma,\tau
(1-\sigma)}X_{\alpha,\tau\sigma}$
satisfy
\[
X^{(\sigma)}_{\alpha,\tau}\stackrel{d}{=}\beta_{\tau\sigma,\tau
(1-\sigma
)}X_{\alpha,\tau\sigma}\stackrel{d}{=}
\beta^{1/\alpha}_{({\tau\sigma}/{\alpha},{(\tau
(1-\sigma
))}/{\alpha})}X_{\alpha,\tau}\stackrel{d}{=}M_{\tau}(F_{X_{\alpha
}\xi
_{\sigma}}).
\]
Which leads to the identity
%
%
\begin{equation}
\label{localtimeid}
\frac{\beta_{\tau\sigma,\tau(1-\sigma)}}{S_{\alpha,\tau\sigma
}}\stackrel
{d}{=}\frac{\beta^{1/\alpha}_{({\tau\sigma}/{\alpha},
({\tau
(1-\sigma)})/{\alpha})}}{S_{\alpha,\tau}}.
\end{equation}
\end{prop}
\begin{pf}
The result is easily checked by following arguments similar
to those
used to establish (\ref{Xid2}). Hence we just note that one uses
the calculation, $\mathcal{C}_{\tau}(\lambda;
X^{(\sigma)}_{\alpha,\tau})=\mathcal{C}_{\tau\sigma}(\lambda;
X_{\alpha,\tau\sigma})$, in place of (\ref{powerid}). The
equality (\ref{localtimeid}) follows immediately since stable
random variables $S_{\alpha}$ are simplifiable (see \cite{Chaumont},
pages 11 and~12).
\end{pf}

Note that Propositions \ref{prop21} and \ref{prop22} show that
\[
-\log\mathcal{C}_{\theta}(\lambda;X_{\alpha,\theta})=\frac
{\theta
}{\alpha}\log(1+\lambda^{\alpha})=\theta
\mathbb{E}[\log(1+\lambda X_{\alpha})].
\]

\subsection{Densities and explicit mixture representations}\label{sec22}

We first describe some more pertinent features of $X_{\alpha}$
(see also \cite{BFRY,JamesBernoulli,JRY}).
\begin{prop}\label{prop23}
Let $X_{\alpha}\stackrel{d}{=}S_{\alpha}/S'_{\alpha}$, having
density (\ref{denX}). Then:
\begin{longlist}
\item the cdf of $X_{\alpha}$ can be represented explicitly as
%
%
\begin{equation}
\label{cdfXY}
F_{X_{\alpha}}(x)=1-\frac{1}{\pi\alpha}\cot^{-1} \biggl(\cot(\pi
\alpha)+\frac{x^{\alpha}}{\sin(\pi\alpha)} \biggr);
\end{equation}
\item its inverse is given by
%
%
\begin{equation}\label{invX}
F^{-1}_{X_{\alpha}}(y)={ \biggl[\frac{\sin(\pi\alpha(y))}{\sin(\pi
\alpha(1-y))} \biggr]}^{1/\alpha};
\end{equation}
\item equations (\ref{cdfXY}) and (\ref{invX}) yield the
identity
%
%
\begin{eqnarray}\label{sinXid}
\sin(\pi\alpha
F_{X_{\alpha}}(y))&=&y^{\alpha}\sin\bigl(\pi
\alpha\bigl(1-F_{X_{\alpha}}(y)\bigr)\bigr)\nonumber\\[-8pt]\\[-8pt]
&=&\frac{y^{\alpha}\sin(\pi\alpha)}{{[y^{2\alpha}+2y^{\alpha
}\cos(\pi
\alpha)+1]}^{1/2}};\nonumber
\end{eqnarray}
\item additionally,
\[
\cos(\pi\alpha
F_{X_{\alpha}}(y))=\frac{y^{\alpha}\cos(\pi\alpha
)+1}{{[y^{2\alpha
}+2y^{\alpha}\cos(\pi
\alpha)+1]}^{1/2}}.
\]
\end{longlist}
\end{prop}
\begin{pf}
This derivation of the cdf is influenced by arguments in
Fujita and Yor \cite{FY} where it becomes clear that it is easier
to work with the density of ${(X_{\alpha})}^{\alpha}$.
Specifically the density of ${(X_{\alpha})}^{\alpha}$ is given by
\[
\frac{\sin(\pi\alpha)}{\pi\alpha}\frac{1}{y^{2}+2y\cos(\pi
\alpha)+1}\qquad\mbox{for }y>0.
\]
It it then easy to obtain the form of the cdf of
${(X_{\alpha})}^{\alpha}$ by direct integration. Now using the
fact that this equates with $F_{X_{\alpha}}(y^{1/\alpha}) $ yields
statement (i). Statement (ii) then follows by using properties
of the inverse cotangent. In order to establish (iii),
use (\ref{invX}) which yields the identity
%
%
\begin{equation}\label{inv}
y=F^{-1}_{X_{\alpha}}(F_{X_{\alpha}}(y))=\biggl[\frac{\sin(\pi
\alpha(F_{X_{\alpha}}(y)))}{\sin(\pi
\alpha(1-F_{X_{\alpha}}(y)))} \biggr]^{1/\alpha}.
\end{equation}
Hence statement (ii) follows.
\end{pf}

We now focus on obtaining explicit distributional formulae for the
pertinent random variables based on their representations as
Dirichlet means. In relation to this, Proposition \ref{prop23} gives precise
details on the pertinent cdf $F_{X_{\alpha}};$ it then remains to
obtain a nice expression for the quantity
\[
\Phi_{\alpha}(x):=\Phi_{X_{\alpha}}(x)=\mathbb{E}[{\log}
|x-X_{\alpha}|]
\]
for $x>0$. The key to calculating $\Phi_{\alpha}(x)$ is the fact
that we showed that $X_{\alpha}$ is a mean functional of the type
$M_{1}(F_{\xi_{\alpha}X_{\alpha}})$, as described in Proposition
\ref{prop21}. This sets up an equivalence between the form of the density
of $X_{\alpha}$ obtained by Lamperti \cite{Lamperti} and that of
$M_{1}(F_{\xi_{\alpha}X_{\alpha}})$, obtained from (\ref{M1}).
Hence we have the following calculation:
\begin{prop}\label{prop24}
For $0<\alpha<1$, and $x>0$,
\[
\Phi_{\alpha}(x)=\frac{1}{2\alpha}\log\bigl(x^{2\alpha}+2x^{\alpha
}\cos
(\alpha
\pi)+1\bigr).
\]
\end{prop}
\begin{pf}
Since $X_{\alpha}\stackrel{d}{=}M_{1}(F_{X_{\alpha}\xi_{\alpha}})$, it
follows by using (\ref{M1}) that the density of $X_{\alpha}$
satisfies the equivalence
\[
f_{X_{\alpha}}(x)=\frac{1}{\pi}\sin\bigl(\pi\alpha
[1-F_{X_{\alpha}}(x)]\bigr)e^{-\alpha\Phi_{\alpha}(x)}x^{\alpha-1}.
\]
Where on the left-hand side we use the expression in (\ref{denX}).
Now applying the identity in (\ref{sinXid}) shows that
\[
f_{X_{\alpha}}(x)=\frac{1}{\pi}\frac{x^{\alpha-1}\sin(\pi\alpha
)}{{[x^{2\alpha}+2x^{\alpha}\cos(\pi
\alpha)+1]}^{1/2}}e^{-\alpha\Phi_{\alpha}(x)}.
\]
Solving this expression for $\Phi_{\alpha}(x)$ concludes the
result.
\end{pf}

Set
%
%
\begin{equation} \label{short}
\rho_{\alpha,\tau}(x^{\alpha})=\frac{\tau}{\alpha}\arctan
\biggl(\frac
{\sin(\pi
\alpha)}{\cos(\pi\alpha)+x^{\alpha}} \biggr)=\pi\tau
[1-F_{X_{\alpha}}(x)]
\end{equation}
and define the function
%
\begin{equation} \label{density}
\Delta_{\alpha,\tau}(x)=\frac{x^{\tau-1}}{\pi}\frac{\sin
(\rho
_{\alpha,\tau}(x^{\alpha}) )}
{{[x^{2\alpha}+2x^{\alpha}\cos(\alpha
\pi)+1]}^{{\tau}/({2\alpha})}}.
\end{equation}

We next obtain density formula for a key class of random variables
that includes the case of $X_{\alpha}$, and $X_{\alpha,1}$.
\begin{prop}\label{prop25}
For $0<\sigma\leq1$, and $x>0$, the
densities of the random variables
\[
X^{(\sigma)}_{\alpha,1}\stackrel{d}{=}\beta_{\sigma,1-\sigma
}X_{\alpha
,\sigma}\stackrel{d}{=}\bigl[\beta_{({\sigma}/{\alpha},
({1-\sigma})/{\alpha})}\bigr]^{{1/\alpha}}X_{\alpha,1}
\]
with $X_{\alpha}\stackrel{d}{=}X^{(\alpha)}_{\alpha,1}$ and
$X_{\alpha,1}\stackrel{d}{=}X^{(1)}_{\alpha,1}$, can be expressed as
$\Delta_{\alpha,\sigma}(x), $ given in (\ref{density}).
Furthermore,
\begin{longlist}
\item
$X^{(\sigma)}_{\alpha,1}\stackrel{d}{=}F^{-1}_{X_{\alpha
}}(U_{\alpha
,\sigma})$
where
$U_{\alpha,\sigma}\stackrel{d}{=}F_{X^{\alpha}_{\alpha
}}([X^{(\sigma
)}_{\alpha,1}]^{\alpha})$
has density
\[
\biggl[\frac{\sin(\pi\alpha)}{\sin(\pi\alpha
u)}\biggr]^{({\alpha-\sigma})/{\alpha}}\frac{\sin(\pi\sigma
(1-u))}{\sin(\pi\alpha(1-u))},\qquad
0<u<1.
\]
\item If $0<\sigma\leq\alpha$, then
$X^{(\sigma)}_{\alpha,1}\stackrel{d}{=}[\beta_{({\sigma
}/{\alpha
},({\alpha-\sigma})/{\alpha})}]^{{1/\alpha}}X_{\alpha}$.
\end{longlist}
\end{prop}
\begin{pf}
The representations of $X^{(\sigma)}_{\alpha,1}$ is
just a special case of Proposition \ref{prop22}. The density is calculated
based on Proposition \ref{prop24} and the results discussed
in \cite{JamesBernoulli} and \cite{CifarelliRegazzini}, as
mentioned previously. Statement (i) takes advantage of the
properties of $F_{X_{\alpha}}$ and is otherwise straightforward to
obtain. Statement (ii) is just a manipulation of the beta random
variables.
\end{pf}

One important aspect of the previous result is that we can use it
to obtain density/mixture representations for the following
Lamperti random variables. This is facilitated by identity
(\ref{Xid2}).
\begin{prop}\label{prop26}
Suppose that $0\leq\theta\leq1-\alpha$, then
$\alpha\leq\sigma^{*}=\theta+\alpha\leq1$ and there is the
distributional identity
\[
X_{\alpha,\theta}\stackrel{d}{=}\beta_{\theta+\alpha,1-\alpha
}X_{\alpha
,\theta+\alpha}\stackrel{d}{=}\beta_{1,\theta}X^{(\sigma
^{*})}_{\alpha,1}.
\]
In particular,
$X_{\alpha,1-\alpha}\stackrel{d}{=}\beta_{1,1-\alpha}X_{\alpha,1}$.
\begin{longlist}
\item
Hence for $0<\theta\leq1-\alpha$, the density of
$X_{\alpha,\theta}$ can be written as
%
%
\begin{equation}\label{shortid1}
f_{X_{\alpha,\theta}}(x)=\theta\int_{0}^{1}
\frac{\Delta_{\alpha,\sigma^{*}}(x/u)}{u{(1-u)}^{1-\theta}}\,du,\qquad
x>0,
\end{equation}
where $\Delta_{\alpha,\sigma^{*}}(x)\ge0$ is the density of
$X^{(\sigma^{*})}_{\alpha,1}$. When $\theta=0$, the density is
$\Delta_{\alpha,\alpha}(x)$ equating with (\ref{denX}).
\item As a special case, when $\alpha\leq1/2$,
$X_{\alpha,\alpha}\stackrel{d}{=}B_{1,\alpha}X^{(2\alpha)}_{\alpha,1}$,
where $X^{(2\alpha)}_{\alpha,1}$ has density
%
%
\begin{equation} \label{shortid3}
\frac{\sin(\pi\alpha)}{\pi}\frac{2\alpha x^{2\alpha
-1}[\cos(\pi\alpha)+x^{\alpha}]}{{[x^{2\alpha}+2x^{\alpha}\cos
(\alpha
\pi)+1]}^{2}}.
\end{equation}
\end{longlist}
\end{prop}
\begin{pf}
The result follows from Propositions \ref{prop22} and \ref{prop25}
by writing
\[
X_{\alpha,\theta}\stackrel{d}{=}\beta_{\theta+\alpha,1-\alpha
}X_{\alpha
,\theta+\alpha}\stackrel{d}{=}\beta_{1,\theta}\beta_{\theta
+\alpha
,1-(\theta+\alpha)}X_{\alpha,\theta+\alpha}.
\]
The simplification in (\ref{shortid3}) follows from
\[
\sin\bigl(2\pi\alpha[1-F_{X_{\alpha}}(x)]\bigr)=\frac{\sin(2\pi\alpha
)+2x^{\alpha
}\sin(\pi
\alpha)}{1+2x^{\alpha}\cos(\pi\alpha)+x^{2\alpha}}.
\]
\upqed\end{pf}

The previous results allow one to obtain simple mixture
representations or densities for $X_{\alpha,\theta}$ in the range
$0\leq\theta\leq1-\alpha$, and $\theta=1$. The fact that we
obtain such results for a continuous range of $\theta$ is
significant, as shown in the next result.
\begin{prop}\label{prop27}
Set $\theta=\sum_{j=1}^{k}\theta_{j}$ where
$\theta_{j}>0$. Furthermore, let $(D_{1},\break\ldots, D_{k})$ denote a
Dirichlet random vector having density proportional to
$\prod_{i=1}^{k}x^{\theta_{i}}_{i}$. That is each
$D_{i}\stackrel{d}{=}\beta_{\theta_{i},\theta-\theta_{i}}$. Then,
\[
X_{\alpha,\theta}\stackrel{d}{=}\sum_{j=1}^{k}D_{j}X_{\alpha
,\theta_{j}},
\]
where $X_{\alpha,\theta_{j}}$ are mutually independent and
independent of $(D_{1},\ldots,D_{k})$. When $\theta_{j}$ are
chosen such that $0<\theta_{j}\leq1-\alpha$, each
$X_{\alpha,\theta_{j}}$ has an explicit density
$f_{X_{\alpha,\theta_{j}}}$ described in (\ref{shortid1}). When
$\theta=k$, one can use $\theta_{j}=1$.
\end{prop}
\begin{pf}
Since we have shown that
$X_{\alpha,\theta}\stackrel{d}{=}M_{\theta}(F_{X_{\alpha}})$, this
result follows directly as a special case of Hjort and
Ongaro (\cite{Hjort}, Proposition 9).
\end{pf}

\subsection{Positive Linnik variables}\label{sec23}

For $\theta>0$,
%
%
\begin{equation}\label{posLinnik}
\chi_{\alpha,\theta}\stackrel{d}{=}\gamma^{1/\alpha}_{{\theta/\alpha}}S_{\alpha}
\end{equation}
denotes the class of generalized Linnik variables as considered in
\cite{Pillai,BondBook,Jayakumar,DevroyeOne,Lin,HT}. The results
in the previous section depend on the validity of the transforms
in (\ref{Cauchy1}) and (\ref{Cauchy2}). It is evident, and known,
that $(1+\lambda^{\alpha})^{-\theta/\alpha}$ appearing in
(\ref{Cauchy1}) is the Laplace transform of $\chi_{\alpha,\theta}$
for $\theta>0$. Hence\vspace*{-1pt} (\ref{Cauchy1}) is verified if one shows
that
$\chi_{\alpha,\theta}\stackrel{d}{=}\gamma_{\theta}X_{\alpha
,\theta}$,
for $\theta>0$. It is already known, using a result of
Devroye \cite{DevroyeLinnik} combined with (\ref{keyid1}), that
%
%
\begin{equation}
\chi_{\alpha,\alpha}\stackrel{d}{=}\gamma^{1/\alpha}_{1}S_{\alpha
}=\gamma_{1}X_{\alpha}\stackrel{d}{=}\gamma_{\alpha}X_{\alpha
,\alpha}.
\end{equation}
Furthermore,\vspace*{-1pt} from Bondesson (\cite{BondBook}, page 38), it follows that
$\chi_{\alpha,\theta}$ are $\operatorname{GGC}(\theta$,\break
$F_{X_{\alpha}})$. In the next result we will verify the usage of
(\ref{Cauchy1}), and use the $X^{(\sigma)}_{\alpha,1}$ to obtain
explicit density representations. In this regard, it is is important to
note that we do not need explicit results for $X_{\alpha,\theta}$ to
get corresponding results for $\chi_{\alpha,\theta}$. In addition we
obtain some interesting identities.
\begin{prop}\label{prop28}
For all $\theta>0$, $\chi_{\alpha,\theta}$ is a
$\operatorname{GGC}(\theta,F_{X_{\alpha}})$ variable that satisfies
\[
\chi_{\alpha,\theta}\stackrel{d}{=}\gamma_{\theta}X_{\alpha
,\theta}.
\]
For $0<\theta=\sigma\leq1$, $
\chi_{\alpha,\sigma}\stackrel{d}{=}\gamma_{1}X^{(\sigma)}_{\alpha,1}
$ and hence has the density
%
%
\begin{equation}\label{extend}
f_{\chi_{\alpha,\sigma}}(x)=\int_{0}^{\infty}e^{x/y}y^{-1}\Delta_{\alpha,\sigma}(y)\,dy.
\end{equation}
See also (\ref{generaldensity}) for $\theta>0$. Additionally:
\begin{longlist}
\item For $\theta>-\alpha$,
%
%
\begin{equation}
\chi_{\alpha,\theta+\alpha}\stackrel{d}{=}\gamma_{\theta+\alpha
}X_{\alpha,\theta+\alpha}\stackrel{d}{=}\gamma_{1+\theta
}X_{\alpha
,\theta}.
\end{equation}
\item Hence, for $\theta>-\alpha$,
%
%
\begin{equation}\label{gammaid}
\gamma^{1/\alpha}_{({\theta+\alpha})/{\alpha}}=\frac{\gamma
_{\theta
+\alpha}}{S_{\alpha,\theta+\alpha}}\stackrel{d}{=}\frac{\gamma
_{1+\theta
}}{S_{\alpha,\theta}}.
\end{equation}
\item For $-\alpha<\theta\leq k$, $k=0,1,2,\ldots,$
\[
\chi_{\alpha,\theta+\alpha}=\gamma_{k+1}X_{\alpha,k}\beta
^{1/\alpha
}_{(({\theta+\alpha})/{\alpha},({k-\theta})/{\alpha})}.
\]
\item For $\theta=\sum_{i=1}^{k}\theta_{i}>0$,
$\chi_{\alpha,\theta}\stackrel{d}{=}\sum_{i=1}^{k}
\chi_{\alpha,\theta_{i}}$, where $\chi_{\alpha,\theta_{i}}$ are
independent.
\end{longlist}
\end{prop}
\begin{pf}
From (\ref{posLinnik}) and using the identity
\[
e^{-{x^{\alpha}}/{s^{\alpha}}}=\mathbb{E}[e^{-{x/s}S_{\alpha}}],
\]
it is easy to see that the density can be expressed as
\begin{eqnarray*}
f_{\chi_{\alpha,\theta}}(x)& \propto&
x^{\theta-1}\int_{0}^{\infty}e^{-{x^{\alpha}}/{s^{\alpha}}}s^{-\theta}f_{\alpha}(s)\,ds\\
& \propto&
x^{\theta-1}\int_{0}^{\infty}\int_{0}^{\infty}e^{-{xv/s}}
(v/s)^{\theta}v^{-\theta}f_{\alpha}(v)f_{\alpha}(s)\,dv\,ds
\end{eqnarray*}
yielding the equivalence with $\gamma_{\theta}X_{\alpha,\theta}$.
The expression in (\ref{extend}) is due to Proposition \ref{prop25}.
Statement (i) follows from (\ref{keyid1}). Statement (ii)
follows by removing $S_{\alpha}$ which is justified since it is a
simplifiable variable. For (iii), apply Proposition~\ref{prop22}(iv)
follows from infinite divisibility.
\end{pf}
\begin{rem} It is not difficult to show that a general expression for
the density of $\chi_{\alpha,\theta}$, for all $\theta>0$, is
obtained by replacing $\Delta_{\alpha,\sigma}$ by
$\Delta_{\alpha,\theta}$ as follows:
%
%
\begin{equation} \label{generaldensity}
f_{\chi_{\alpha,\theta}}(z)=\frac{1}{\pi}\int_{0}^{\infty}\frac
{e^{-z x}\sin(\pi\theta
F_{X_{\alpha}}(x))\,dx}{{[x^{2\alpha}+2x^{\alpha}\cos(\alpha
\pi)+1]}^{{\theta}/({2\alpha})}}.
\end{equation}
However, $\Delta_{\alpha,\theta}$ can take negative values when
$\theta>1$, so this does not in general yield a mixture
representation for $\chi_{\alpha,\theta}. $ Nonetheless, it may
not be difficult to evaluate numerically which is relevant for
the Mittag--Leffler functions discussed in Section \ref{sec3}.
\end{rem}
\begin{rem}
The gamma identity in statement (\ref{gammaid}) of the previous
proposition is quite remarkable, and, as we shall see below, has
some interesting implications. We note that although not obvious,
our result coincides with a variation of Bertoin and
Yor (\cite{BerYor}, Lemma 6). Checking moments one can see that,
in their notation,
$J_{s,s/\alpha}\stackrel{d}{=}S^{-\alpha}_{\alpha,s}$ for $s>0$ and
for $\theta>-\alpha$, $
J^{(\alpha)}_{1+\theta,(\theta+/\alpha)/\alpha}\stackrel
{d}{=}S^{-\alpha
}_{\alpha,\theta}.
$ Our work provides some additional interpretation of their
variables (see also \cite{JamesYor}). See
Kotlarski \cite{Kotlarski} for a general characterization of cases
where products of variables result in gamma variables.
\end{rem}

\subsection{Generalized Pareto laws}\label{sec24}

Influenced in part by the gamma identity (\ref{gammaid}), we next
look at relationships between the Lamperti laws and a class of
generalized Pareto distributions. We note that the next result
also plays an important role in Section \ref{sec4} when discussing
continuous state branching processes. Define random variables
\[
W^{1/\alpha}_{\alpha,\theta}:=\biggl(\frac{U^{{\alpha
}/{\theta
}}}{1-U^{{\alpha/\theta}}} \biggr)^{{1/\alpha}}.
\]
These represent a sub-class of generalized Pareto distributions
with cdf and density given as
\[
F_{W^{1/\alpha}_{\alpha,\theta}}(y)=\frac{y^{\theta}}{(1+y^{\alpha
})^{\theta/\alpha}};\qquad
f_{W^{1/\alpha}_{\alpha,\theta}}(y)=\frac{\theta
y^{\theta-1}}{{(1+y^{\alpha})}^{(\theta+\alpha)/\alpha}}.
\]
\begin{prop}\label{prop29}
Let $U$ denote a Uniform $[0,1]$ random
variable, then
for $\theta>0$:
\begin{longlist}
\item There is the identity
\[
W^{1/\alpha}_{\alpha,\theta}:= \biggl(\frac{U^{{\alpha
}/{\theta
}}}{1-U^{{\alpha}/{\theta}}} \biggr)^{{1}/{\alpha}}
\stackrel{d}{=}{ \biggl(\frac{\gamma_{{\theta}/{\alpha
}}}{\gamma
_{1}} \biggr)}^{{1}/{\alpha}}
\stackrel{d}{=}\frac{\chi_{\alpha,\theta}}{\gamma_{1}}
\stackrel{d}{=}\frac{\gamma_{\theta}}{\gamma_{1}}X_{\alpha,\theta}.
\]
\item For $0<\sigma\leq1$, the random variable
$\Sigma_{\alpha,\sigma}\stackrel{d}{=}\gamma_{1}/X^{(\sigma
)}_{\alpha,1}$
has Laplace transform
%
%
\begin{equation}\label{CSBPhook}
\mathbb{E}[e^{-\lambda
\Sigma_{\alpha,\sigma}}]=1-\lambda^{\sigma}{(1+\lambda^{\alpha
})}^{-\sigma/\alpha}.
\end{equation}
\end{longlist}
\end{prop}
\begin{pf}
Statement (i) is an application of (\ref{gammaid}).
For statement (ii) notice that
\[
\mathbb{P}\biggl(\frac{\gamma_{1}}{\Sigma_{\alpha,\sigma}}>\lambda
\biggr)=\mathbb
{E}[e^{-\lambda\Sigma_{\alpha,\sigma}}],
\]
but this is the survival function of the random variable
\[
\frac{\gamma'_{1}}{\gamma_{1}}X^{(\sigma)}_{\alpha,1}\stackrel
{d}{=}\frac
{\gamma_{\sigma}}{S_{\alpha,\sigma}}\frac{S_{\alpha}}{\gamma
_{1}}\stackrel{d}{=}W^{1/\alpha}_{\alpha,\sigma}.
\]
\upqed\end{pf}
\begin{rem} As we shall discuss in Section \ref{sec4}, Statement (ii), (\ref
{CSBPhook}) serves to identify
explicitly the (unknown) limiting distribution obtained by
\cite{ZolotarevSlack} and \cite{Slack} corresponding to
$\sigma=1$. It is relevant also to note that the density of
$W^{1/\alpha}_{\alpha,1}$ is the only case that corresponds to a
Laplace transform. So here we see a distinguishing feature of
$X_{\alpha,1}$.
\end{rem}

\subsection{$z$-variables and hyperbolic laws}\label{sec25}

Proposition \ref{prop29}, along with the works of \cite
{Biane,PYhyper,YanoYor,JeanHu}, motivate us to consider several
questions related to $z$-distributions which are distributed as
the logarithm of of the ratio of independent gamma variables. We
also believe that some of the variables we discuss will be of
interest in terms of applications along the lines discussed
in \cite{BNKent} and \cite{JeanHu}. In fact, \cite{JeanHu}
suggests the use of a class of variables that turn out to be
equivalent in distribution to $\log(X_{\alpha})$. Naturally we do
this in the spirit of highlighting what one can do with Lamperti
laws. We also obtain additional information about these variables.

In particular, for illustration, we consider the following generic
type of problem. Suppose for generic variables $X,Y,Z$ with $Z$
and $Y$ independent there is the following relation:
\[
X\stackrel{d}{=}Y+Z.
\]
One natural question is to ask, given explicit information about
$X$ and $Y$, what $Z$ satisfies the above equation? In addition,
does $Z$ have an explicit density or mixture representation?
Notice also that if the characteristic function of $Z$ is not
known then we can use $X$ and $Y$ to obtain this. We will consider
$Z$ that are variants of Lamperti laws.

We first give a brief discussion on $z$-distributions.
Following \cite{YanoYor}, the class of $z$-distributions are
defined as
$\pi^{-1}\log(\gamma_{\theta_{1}}/\gamma_{\theta_{2}})$, having
characteristic function
\[
\mathbb{E}\bigl[e^{{i\lambda}/{\pi}
\log(\gamma_{\theta_{1}}/\gamma_{\theta_{2}})}\bigr]=\frac{\Gamma
(\theta
_{1}+i{\lambda}/{\pi})\Gamma(\theta_{2}-{i\lambda}/{\pi})}
{\Gamma(\theta_{1})\Gamma(\theta_{2})}.
\]
As special cases, the variables, for $0<\sigma<1$,
$M_{\sigma}\stackrel{d}{=}\pi^{-1}\log(\gamma_{\sigma}/\gamma
_{1-\sigma})$,
have Meixner distributions with characteristic function
%
%
\begin{equation} \label{Meixnercf}
\mathbb{E}[e^{-i\lambda
M_{\sigma}}]=\frac{\cos(\varepsilon_{\sigma})} {\cosh
(\lambda-i\varepsilon_{\sigma})},
\end{equation}
where $\varepsilon_{\sigma}=\pi(\sigma-1/2)$. Note that a Meixner
distributed random variable is usually defined as
$(1/2)M_{\sigma}$.
$\mathbb{S}_{1}\stackrel{d}{=}\pi^{-1}\log(U/(1-U))$ has a \textit{logistic
distribution} with the characteristic function
\[
\mathbb{E}[e^{i\lambda
\mathbb{S}_{1}}]=\frac{\lambda}{\sinh(\lambda)}
\]
and
$\mathbb{C}_{1}\stackrel{d}{=}\pi^{-1}\log(\gamma'_{{1/2}}/\gamma_{{1/2}})$
has the \textit{hyperbolic distribution} with characteristic
function
\[
\mathbb{E}[e^{i\lambda
\mathbb{C}_{1}}]=\frac{1}{\cosh(\lambda)}.
\]

It is known that the variables $\mathbb{S}_{1}$ and
$\mathbb{C}_{1}$ satisfy
%
%
\begin{equation}\label{sumC}
\mathbb{C}_{1}\stackrel{d}{=}\mathbb{S}_{1}+\mathbb{T}_{1},
\end{equation}
where $\mathbb{T}_{1}$ is an independent variable having characteristic
function
\[
\mathbb{E}[e^{i\lambda
\mathbb{T}_{1}}]=\frac{\tanh(\lambda)}{\lambda}.
\]
Biane and Yor \cite{BianeYor} showed that the density of
$\mathbb{T}_{1}$ is
\[
f_{\mathbb{T}_{1}}(x)=\frac{1}{\pi}\log
\biggl(\coth\biggl(\frac{4}{\pi}|x|\biggr)\biggr),\qquad
-\infty< x<\infty.
\]
\begin{rem}
Note that the characteristic function of $\alpha
\pi^{-1}\log(X_{\alpha})$ is equivalent to
\[
\mathbb{E}\bigl[e^{ {i
\alpha\lambda}/{\pi}\log(X_{\alpha})}\bigr]=
\frac{\sinh(\alpha\lambda)}{\alpha\sinh(\lambda)}.
\]
This
expression can be found in Chaumont and Yor (\cite{Chaumont}, page
147). In addition we see that this characteristic function agrees
with the class of generalized secant hyperbolic distributions
discussed, for instance, in \cite{JeanHu}. See also \cite{PYCauchy}
for more on the variables $\mathbb{C}_{1},\mathbb{T}_{1}$ and
$\mathbb{S}_{1}$.
\end{rem}

In the next result we obtain a description of the characteristic
function of $\log(X_{\alpha,\theta})$ and
$\log(X^{(\sigma)}_{\alpha,1})$.
\begin{prop}\label{prop210}
From Proposition \ref{prop29}, it follows that:
\begin{longlist}
\item For $\theta>0$,
\[
\frac{1}{\alpha}\log \biggl(\frac{\gamma_{{\theta}/{\alpha
}}}{\gamma
_{1}} \biggr)
\stackrel{d}{=}
\log\biggl(\frac{\gamma_{\theta}}{\gamma_{1}} \biggr)+\log
(X_{\alpha
,\theta}).
\]
\item For $\theta>-\alpha$
\[
\mathbb{E}\bigl[e^{i\lambda
\pi^{-1}\log(X_{\alpha,\theta})}\bigr]=\frac{\Gamma(({\theta
+\alpha
})/{\alpha}+i{\lambda}/({\alpha\pi}))\Gamma(1-{i\lambda
}/({\alpha\pi
}))\Gamma(1+\theta)}
{\Gamma(1+\theta+i{\lambda}/{\pi})\Gamma(1-{i\lambda
}/{\pi})
\Gamma(({\theta+\alpha})/{\alpha})}.
\]
\item For $0<\sigma\leq1$,
\[
\frac{1}{\alpha}\log \biggl(\frac{\gamma_{{\sigma}/{\alpha
}}}{\gamma
_{1}} \biggr)
\stackrel{d}{=}
\log\biggl(\frac{\gamma'_{1}}{\gamma_{1}} \biggr)+\log\bigl(X^{(\sigma
)}_{\alpha,1}\bigr),
\]
where $\log(X^{(\sigma)}_{\alpha,1}) $ has density
\[
\frac{1}{\pi}\frac{\sin(\rho_{\alpha,\sigma}(e^{z\alpha}) )}
{{[e^{-2z\alpha}+2e^{-z\alpha}
\cos(\alpha\pi)+1]}^{{\sigma}/({2\alpha})}},\qquad
-\infty<z<\infty,
\]
and characteristic function
\[
\mathbb{E}\bigl[e^{i\lambda\log(X^{(\sigma)}_{\alpha,1})}\bigr]=\frac{\sinh(\lambda
\pi
)}{\lambda\pi}\frac{\Gamma({\sigma}/{\alpha}+{i\lambda
}/{\alpha
})\Gamma(1-
{i\lambda}/{\alpha})}{\Gamma({\sigma}/{\alpha})}.
\]
\end{longlist}
\end{prop}
\begin{pf}
This follows as a simple consequence of our previous
results and the characteristic functions of $z$-distributions.
\end{pf}

The next result identifies some variables that have characteristic
functions based on hyperbolic functions and also have explicit
densities. Define
\[
H_{\alpha,\sigma}\stackrel{d}{=}\frac{X^{(\alpha\sigma)}_{\alpha
,1}}{\beta^{1/\alpha}_{1-\sigma,\sigma}}\stackrel{d}{=}
\frac{\beta^{1/\alpha}_{\sigma,1-\sigma}}
{\beta^{1/\alpha}_{1-\sigma,\sigma}}X_{\alpha},
\]
where the equality follows from (\ref{localtimeid}). In addition
for $\alpha\delta\leq\theta\leq\alpha(1-\delta)$, for
$\delta\leq1/2$, define
\[
L^{(\delta)}_{\alpha,\theta}\stackrel{d}{=} \biggl(\frac{\beta
_{(\delta
,({\theta-\alpha\delta})/{\alpha})}}
{\beta_{(({1-\alpha(1-\delta)})/{\alpha},({\alpha(1-\delta
)-\theta
})/{\alpha})}} \biggr)^{{1/\alpha}}
\frac{S_{\alpha,1-\theta}}{S_{\alpha,\theta}},
\]
where one can easily check that the density of
$S_{\alpha,1-\theta}/S_{\alpha,\theta}$, denoted as
$f_{1-\theta,\theta}$, satisfies
\[
f_{1-\theta,\theta}(x)=
c_{\alpha,\theta}x^{-(1-\theta)}f_{X_{\alpha,1}}(x)=c_{\alpha
,\theta
}x^{-(1-\theta)}\Delta_{\alpha,1}(x)
\]
for
\[
c_{\alpha,\theta}=\frac{\Gamma(1/\alpha+1)\Gamma(\theta+1)\Gamma
(2-\theta)}{\Gamma(\theta/\alpha+1)\Gamma((1-\theta)/\alpha+1)}.
\]
\begin{prop}\label{prop211} For $0<\sigma<1$, $\delta\leq1/2$ and
$\alpha\delta\leq\theta\leq\alpha(1-\delta)$, there are the
following relationships:
\begin{longlist}
\item
$
\frac{1}{\alpha}\log{ (\frac{\gamma_{\sigma}}{\gamma
_{1-\sigma
}} )}
\stackrel{d}{=}
\log(\frac{\gamma'_{1}}{\gamma_{1}} )+\log(H_{\alpha
,\sigma}).
$ Hence,
\[
\mathbb{E}\bigl[e^{{i\alpha\lambda}/{\pi}\log(H_{\alpha,\sigma})}\bigr]=\frac
{\cos
(\varepsilon_{\sigma})\sinh(\lambda\alpha)}
{\lambda\alpha\cosh(\lambda-i\varepsilon_{\sigma})}
=\frac{\sinh(\alpha\lambda)}{\alpha\sinh(\lambda)}
\frac{\cos(\varepsilon_{\sigma})\sinh(\lambda)} {\lambda\cosh
(\lambda-i\varepsilon_{\sigma})},
\]
where $\varepsilon_{\sigma}=\pi(\sigma-1/2)$.
\item
$
\frac{1}{\alpha}\log{ (\frac{\gamma_{\delta}}{\gamma
_{1-\delta
}} )}
\stackrel{d}{=}
\log(\frac{\gamma_{\theta}}{\gamma_{1-\theta}} )+\log
(L^{(\delta)}_{\alpha,\theta}).
$ With
\[
\mathbb{E}\bigl[e^{{i\alpha\lambda}/{\pi}\log(L^{(\delta)}_{\alpha,\theta
})}\bigr]=\frac
{\cos(\varepsilon_{\delta})\cosh
(\lambda\alpha-i\varepsilon_{\theta})}{\cos(\varepsilon_{\theta
})\cosh
(\lambda-i\varepsilon_{\delta})}.
\]
\item When $\delta=1/2$, then $\theta=\alpha/2$, and
$L^{(1/2)}_{\alpha,\alpha/2}\stackrel{d}{=}S_{\alpha,1-\alpha
/2}/S_{\alpha,\alpha/2}$
where $\log(L^{(1/2)}_{\alpha,\alpha/2})$, has density
\[
\frac{c_{\alpha,\alpha/2}}{\pi}\frac{e^{z{\alpha}/{2}}\sin(\rho_{\alpha,1}
(e^{z\alpha}) )}
{{[e^{2z\alpha}+2e^{z\alpha}\cos(\alpha\pi)+1]}^{{1/(2\alpha)}}},\qquad
-\infty<z<\infty,
\]
and characteristic function
\[
\mathbb{E}\bigl[e^{{i\lambda}/{\pi}\log(L^{(1/2)}_{\alpha,\alpha/2})}\bigr]=\frac
{\cosh
(\lambda-i\varepsilon_{\alpha/2})}{\cos(\varepsilon_{\alpha
/2})\cosh
(\lambda/\alpha)}.
\]
\end{longlist}
\end{prop}
\begin{pf}
All the characteristic functions follow from that of
$\mathbb{S}_{1}$ and the Meixner distributions (\ref{Meixnercf}).
In order to establish (i) we use the identity
$\gamma_{1}/S_{\alpha}=\gamma^{1/\alpha}_{1}$ and apply this to
the second definition of $H_{\alpha,\sigma}$. Statement (ii)
follows from a manipulation of (\ref{gammaid}) to force the form
of the first two variables. The choice of $1-\theta$ and $\theta$
in the definition of $L^{(\delta)}_{\alpha,\theta}$ was
deliberately made so that we could get an explicit expression of
the density of the relevant ratios of stable variables.
\end{pf}

\section{Mittag--Leffler functions}\label{sec3}

In this section we obtain integral
representations and other identities
for a generalization of the Mittag--Leffler function given by
%
%
\begin{equation}\label{altM}
\mathrm{E}^{(\theta/\alpha+1)}_{\alpha,1+\theta}(-\lambda)=\sum
_{k=0}^{\infty}\frac{{(-\lambda)}^{k}}{k!}\frac{{(\theta/\alpha
+1)}_{k}}{\Gamma(\alpha
k+\theta+1)}\qquad\mbox{for }\theta>-\alpha,
\end{equation}
where
\[
(\theta/\alpha+1)_{k}=\frac{\Gamma(\theta/\alpha+1+k)}{\Gamma
(\theta
/\alpha+1)}.
\]
So when $\theta=0$, one recovers the Mittag--Leffler function as
\[
\mathrm{E}_{\alpha,1}(-\lambda)=\mathrm{E}^{(1)}_{\alpha,1}(-\lambda
)=\mathrm{E}^{(0)}_{\alpha,0}(-\lambda).
\]

Note that using simple cancelations involving gamma functions it
is easy to show that for $\theta>0$,
%
%
\begin{equation}\label{altM2}
\mathrm{E}^{(\theta/\alpha+1)}_{\alpha,1+\theta}(-z)=\frac{\Gamma
(\theta
)}{\Gamma(1+\theta)}\mathrm{E}^{(\theta/\alpha)}_{\alpha,\theta}(-z).
\end{equation}

Recall that the Mittag--Leffler function can be expressed as the
Laplace transform of $S^{-\alpha}_{\alpha}$. One can show that by
taking a Taylor expansion and calculating moments that, for
$\theta>-\alpha$,
%
%
\begin{equation}\label{MittagLeffler}
\mathbb{E}[e^{-z/S^{\alpha}_{\alpha,\theta}}]
=\mathbb{E}[e^{-z^{1/\alpha}X_{\alpha,\theta}}]
=\Gamma(1+\theta)\mathrm{E}^{(\theta/\alpha+1)}_{\alpha,1+\theta}(-z).
\end{equation}
One consequence of this observation is that one can use a
Monte Carlo method based on $S_{\alpha,\theta}$ to evaluate this
quantity. The next result develops more connections with
$X_{\alpha,\theta}$ and $\chi_{\alpha,\theta}$.
\begin{prop}\label{prop31} For $\theta>-\alpha$:
\begin{longlist}
\item $f_{1/X_{\alpha,\theta}}(x)=x^{-\theta}f_{X_{\alpha
,\theta}}(x);$
\item $f_{\gamma_{1}/X_{\alpha,\theta}}(x)=\Gamma(\theta
+1)x^{-\theta}f_{\chi_{\alpha,\theta+\alpha}}(x)$;
\item
for $\theta>0$,
\[
\mathbb{P}\biggl(\frac{\gamma_{1}}{X_{\alpha,\theta}}>x\biggr)
=\mathbb{E}[e^{-xX_{\alpha,\theta}}]
=\Gamma(\theta)\mathrm{E}^{(\theta/\alpha)}_{\alpha,\theta}(-x^{\alpha
})=\Gamma(\theta)x^{1-\theta}f_{\chi_{\alpha,\theta}}(x).
\]
Hence if $\theta=\sigma$ then these expressions are explicitly
determined by (\ref{extend}), otherwise one might
use (\ref{generaldensity}).
\item Applying Proposition \ref{prop27}, it follows that for $\sum
_{i=1}^{k}\theta_{i}=\theta$, where $\theta_{i}>0$,
\[
\mathrm{E}^{(\theta/\alpha)}_{\alpha,\theta}(-z)=\int_{\mathcal
{S}_{k}}\prod_{i=1}^{k}\mathrm{E}^{(\theta_{i}/\alpha)}_{\alpha
,\theta_{i}}(-z
x^{\alpha}_{i})x_{i}^{\theta_{i}-1}\,dx_{i},
\]
where
$\mathcal{S}_{k}=\{(x_{1},\ldots,x_{k})\dvtx 0<\sum_{i=1}^{k}x_{i}\leq
1\}$.
\end{longlist}
\end{prop}

The result is fairly straightforward using the basic definitions
as ratios of stable variables and the identities in Proposition
\ref{prop25}. We omit the details.
\begin{rem} Recall that for $\alpha=1/2$,
$X_{1/2,\theta}\stackrel{d}{=}\gamma_{\theta+1/2}/\gamma_{1/2}$.
Using the notation in Pitman \cite{Pit02}, equations (88), (98) and Lemma
15, and noting \cite{Pit02}, equation (104), the conditional
moments of the meander length $(1-G_{1})$, of Brownian motion on
$[0,1]$ conditioned on its local time $L_{1}$ is related to the
generalized Mittag--Leffler function as follows:
\begin{eqnarray*}
\mathbb{E}_{1/2,0}\bigl[{(1-G_{1})}^{\theta+1/2}|L_{1}
=\sqrt{2}\lambda\bigr]
&=&\mathbb{E}[e^{-({\lambda^{2}}/{2})X_{1/2,\theta}}]\\
&=&\mathbb{E}[{|B_{1}|}^{\theta+1/2}]h_{-(2\theta+1)}(\lambda)\\
&=& \Gamma(1+\theta)\mathrm{E}^{(2\theta+1)}_{1/2,\theta+1}\biggl(-\frac{\lambda}{\sqrt{2}}\biggr)\\
&=&\mu(\theta+1/2\|\lambda),\qquad\mbox{for }\theta>-1/2,
\end{eqnarray*}
which corresponds to the moments of the structural distribution of
the Brownian excursion partition. $h_{-2q}(\lambda)$ is a Hermite
function, here $q=\theta+1/2$.
\end{rem}
\begin{rem} Equation (\ref{altM}) is a special case of the function
introduced by Prabhakar \cite{Prab},
%
%
\begin{equation}\label{Mittaggen2}
\mathrm{E}^{\gamma}_{\rho,\mu}(-\lambda)=\sum_{k=0}^{\infty
}\frac
{{(-\lambda)}^{k}}{k!}\frac{(\gamma)_{k}}{\Gamma(\rho
k+\mu)},
\end{equation}
where $(\rho,\mu,\gamma\in\mathbb{C}, \operatorname{Re}(\rho)>0)$.
Equation (\ref{altM}) is the case where $\gamma=(\theta+\alpha)/\alpha$ and
$\mu=\theta+1$. Additionally, quantity (\ref{altM}) represents
a special sub-class of yet more general Mittag--Leffler-type
functions which are discussed, for instance, in Kilbas, Saigo and
Megumi \cite{KilbasSaigo}. See also \cite
{Anh,BerMittag,BNL,CT,Djr,Hilfer,Kirbook,kir1,KypScale}.
\end{rem}

\section{\texorpdfstring{The explicit L\'evy density of Stable CSBPs and the
Zolotarev--Slack distribution}{The explicit Levy density of Stable CSBPs and the
Zolotarev--Slack distribution}}\label{sec4}

We now show how our results lead to an explicit identification of
the L\'evy density of the semigroup of stable continuous state
branching processes of index $1<\delta<2$, that is,
$\delta=1+\alpha$ and the limiting distributions first obtained by
Zolotarev \cite{ZolotarevSlack} and Slack \cite{Slack}. We also
mention briefly its connection to the work of \cite{BBS,BerLegall} on
beta coalescents. From Lamperti, \cite{LampertiBerk}
continuous state branching processes (CSBP) are Markov processes
that can be characterized as limits of Galton--Watson branching
processes when the population size grows to infinity.
A~$(1+\alpha)$-stable (CSBP) process $(Y_{t}, t>0)$ is a Markov
process whose semigroup is specified by
%
%
\begin{equation}\label{CSBPtransform}
\mathbb{E}_{a}[e^{-\lambda Y_{t}}]=\mathbb{E}[e^{-\lambda
Y_{t}}|Y_{0}=a]=e^{-a\upsilon_{t}(\lambda)},
\end{equation}
where
\[
\upsilon_{t}=\int_{0}^{\infty}(1-e^{-s\lambda})\nu_{t}(ds)=\lambda{(\alpha
t+\lambda^{\alpha})}^{-1/\alpha}.
\]
Furthermore (see, for instance, \cite{BBS,BerLegall}), there exists
a process $(Y(t,a); t>0, a>0)$ such that for each $t$,
$Y(t,\cdot)$ is a compound Poisson process with intensity
$\nu_{t}$. Related to this result, Kawazu and
Watanabe \cite{Kawazu} show that all continuous state branching
processes with immigration arise as limits of Galton--Watson
processes with immigration. Analogous to (\ref{CSBPtransform}),
Theorem 2.3 of their work yields the limiting $(1+\alpha)$-stable
(CSBP) with immigration $(\hat{Y}(t), t>0)$ satisfying
%
%
\begin{equation} \label{CSBPI}
\mathbb{E}_{a}[e^{-\hat{Y}_{t}}]=(1+\alpha\lambda^{\alpha}t)^{-{d}/({\alpha
c})}\mathbb{E}_{a}[e^{-Y_{t}}].
\end{equation}
It is evident from (\ref{CSBPI}) that the entrance laws of
${\hat Y}$ are positive Linnik distributions. However, the
intensity ${\nu_{t}}$, which plays a fundamental role
in \cite{BBS,BerLegall}, is only known up to its Laplace
transform. It is known that this Laplace transform coincides, up
to some scaling factors, with the Laplace transforms of the
limiting distributions obtained by Zolotarev \cite{ZolotarevSlack}
and Slack \cite{Slack}. Before identifying these expressions we
recount the limiting distributions obtained by Slack \cite{Slack}
and Berestycki, Berestycki and Schweinsberg \cite{BBS}. As
described, for instance, in \cite{BBS}, Slack's result describes
the limiting distribution, say $\mu_{\alpha}$, of the number of
offspring in generation $n$ of a critical Galton--Watson process,
re-scaled to have mean~$1$ and conditioned to be positive, when
the offspring distribution is in the domain of attraction of a
stable law of index $1<\delta<2$. This result complements
Yaglom's \cite{Yag} well-known result for the case where the
offspring distribution has finite variance. In that case the
limiting distribution is exponential with mean $1$. Precisely,
following the exposition in \cite{Nag}, we state a variation of
Slack's result.
\begin{prop}[(Slack (1968) \cite{Slack})]\label{prop41} Let $Z=(Z_{n},n>0)$ denote
a super-critical Galton--Watson process initiated by a single
process. Furthermore, suppose the nonextinction probability
$Q_{n}=\mathbb{P}(Z_{n}>0)$, satisfies
\[
Q_{n}=n^{-1/\alpha}L(n),
\]
where $L(x)$ is a slowly varying function. Then,
%
%
\begin{equation} \label{Slack1}
\lim_{n\rightarrow\infty}\mathbb{P}(Q_{n}Z_{n}\leq
x|Z_{n}>0)=\mu_{\alpha}([0,x]),
\end{equation}
where for each $0<\alpha<1$, $\mu_{\alpha}$ is the distribution of
a random variable $\Sigma_{\alpha}$ satisfying
%
%
\begin{equation}\label{Slack2} \int_{0}^{\infty}e^{-\lambda
w}\mu_{\alpha}(dw)=\mathbb{E}[e^{-\lambda
\Sigma_{\alpha}}]=1-\lambda{(1+\lambda^{\alpha})}^{-1/\alpha}.
\end{equation}
\end{prop}

Zolotarev (\cite{ZolotarevSlack}, Theorem 7) also obtained this
limit in the case of a class of continuous parameter regular
branching processes. However, prior to our work, an explicit
description of its density or corresponding random variable was
not known. It is now evident from (\ref{Slack2}) that
$\Sigma_{\alpha}\stackrel{d}{=}\Sigma_{\alpha,1}$ in Proposition \ref{prop29},
as we mentioned previously. These limits and or discussions
related to (\ref{CSBPtransform}), (\ref{CSBPI}) appear more
recently in, for instance, \cite
{BBS,BerLegall,EvansRalph,KypCSBP,Lageras,Patie}. Before we summarize
our results we shall
say a bit more about the context of \cite{BBS}. Random variables
with law $\mu_{\alpha}$ arise in the work of Berestycki,
Berestycki and Schweinsberg \cite{BBS} (see also \cite{BerLegall})
in connection with $\operatorname{Beta}(2-\delta,\delta)$ coalescents for
$1<\delta<2$. See in particular (\cite{BBS}, Theorem 1.2).
Equivalently these are Beta$(1-\alpha,1+\alpha)$ coalescents.

In addition, there is a related result of Berestycki, Berestycki
and Schweinsberg \cite{BBS} that, with some work, would have
otherwise allowed us describe the law $\mu_{\alpha}$. We quote
their result below.
\begin{prop}[(Berestycki, Berestycki and Schweinsberg \cite{BBS}, Proposition~1.5)]\label{prop42}
Let $(\Pi(t),t>0)$ denote a $\operatorname{Beta}(1-\alpha,1+\alpha)$
coalescent where $0<\break\alpha<1$, and let
$K(t)$ denote the asymptotic frequency of the block of $\Pi(t)$
containing $1$. Then
%
%
\begin{equation} \label{Bias}
\bigl(\Gamma(\alpha+2)t^{-1}\bigr)^{{1/\alpha}}K(t)\mathop
{\rightarrow}_{d}
\zeta_{\alpha}\qquad\mbox{as }t\downarrow0,
\end{equation}
where $\zeta_{\alpha}$ is a random variable satisfying
%
%
\begin{equation}\label{Biasd}
\mathbb{E}[e^{-\lambda
\zeta_{\alpha}}]={(1+\lambda^{\alpha})}^{-(1+\alpha)/\alpha}.
\end{equation}
\end{prop}

Furthermore, as noted in \cite{BBS}, $\zeta_{\alpha}$ has the
size biased distribution
%
%
\begin{equation}\label{sizebias}
\mathbb{P}(\zeta_{\alpha}\in dx)=x\mu_{\alpha}(dx).
\end{equation}
We now summarize our result which again demonstrates the relevance
of the random variable
$X_{\alpha,1}\stackrel{d}{=}S_{\alpha}/S_{\alpha,1}$, which has
explicit density
%
%
\begin{equation} \label{density1}
\Delta_{\alpha,1}(x)=\frac{1}{\pi}\frac{\sin(
{1/\alpha
}\arctan({\sin(\pi
\alpha)}/({\cos(\pi\alpha)+x^{\alpha}})))}
{{[x^{2\alpha}+2x^{\alpha}\cos(\alpha
\pi)+1]}^{{1/(2\alpha)}}}.
\end{equation}
\begin{prop}\label{prop43} For $0<\alpha<1$, and
$\Delta_{\alpha,1}(x)$ defined in (\ref{density1}), there are
the following results.
\begin{longlist}
\item The random variable described
in (\ref{Bias}) and (\ref{Biasd}), $\zeta_{\alpha}$, satisfies
\[
\zeta_{\alpha}\stackrel{d}{=}\gamma_{2}X_{\alpha,1}=\chi_{\alpha
,1+\alpha}.
\]
\item Let $\Sigma_{\alpha}$ and $\mu_{\alpha}$ be as in (\ref
{Slack1})
and (\ref{Slack2}), then
\[
\Sigma_{\alpha}\stackrel{d}{=}\frac{\gamma_{1}}{X_{\alpha
,1}}\stackrel
{d}{=}\Sigma_{\alpha,1}.
\]
\item Furthermore, for each $x>0$,
\[
\mathbb{P}(\Sigma_{\alpha}>
x)=\mu_{\alpha}([x,\infty))=\mathrm{E}^{(1/\alpha)}_{\alpha
,1}(-x^{\alpha})=f_{\gamma_{1}X_{\alpha,1}}(x).
\]
\item The
L\'evy density $\nu_{t}$ corresponding to the $(1+\alpha)$ (CSBP)
specified by (\ref{CSBPtransform}) is
\begin{eqnarray*}
\nu_{t}(x)&=&{{(\alpha t)}^{-{(1+\alpha)}/{\alpha}}x^{\alpha-1}}
\mathrm{E}^{({1+\alpha})/{\alpha}}_{\alpha,1+\alpha}\biggl(-\frac{x^{\alpha}}{\alpha t}\biggr)\\
&=& {(\alpha t)}^{-{2}/{\alpha}}\int_{0}^{\infty}e^{-{xy}/{ {(\alpha t)}^{1/\alpha}}}y\Delta_{\alpha,1}(y)\,dy.
\end{eqnarray*}
\end{longlist}
\end{prop}
\begin{pf}
Statements (i) and (ii) are now quite obvious from
Propositions \ref{prop28} and \ref{prop29}. Statement (iii) is deduced from
Propositions \ref{prop25} and \ref{prop31}. Statement (iv) follows from these
facts and the specifications for $\nu_{t}$, described in
\cite{BBS}, Lem\-ma~2.2, that is,
\[
\nu_{t}(x)=(\alpha t)^{-{1/\alpha}}\mu_{\alpha,t}(x),
\]
where $\mu_{\alpha,t}(x)$ is the density of the random variable
$(\alpha t)^{{1/\alpha}}\Sigma_{\alpha}$.
\end{pf}

\section{Occupation times of generalized Bessel bridges}\label{sec5}

We now show how our results for $X_{\alpha,\theta}$ can be used to
obtain new results related to $A^{+}_{1}$, which is equivalent in
distribution to $P_{\alpha,\theta}(p)$, under
$\mathbb{P}^{(p)}_{\alpha,\theta}$. This can be seen as a
continuation of a subset of the work of James, Lijoi and
Pr\"{u}nster \cite{JLP}, who looked at
more general $\operatorname{PD}(\alpha,\theta)$ mean functionals,
where, with the exception of $\alpha=1/2$, the best results for
describing the density of $P_{\alpha,\theta}(p)$ were obtained for
$\theta=1$, and $\theta=1-\alpha$. The results for $\alpha=1/2$
are classic. For $\alpha=1/2$, and $p=1/2$, L\'evy \cite{Levy}
showed that $A^{+}_{1}$ under $(1/2,0)$ and $(1/2,1/2)$, follow
the Arcsine and $\operatorname{Uniform}[0,1]$ distributions, respectively. A
general formula for $(1/2,\theta)$, for all $\theta>-1/2$ can be
found in Carlton \cite{Carlton}, equation (3.4) (see
also \cite{Keilson}) and is given by
\[
\mathbb{P}^{(p)}_{1/2,\theta}(A^{+}_{1}\in
dy)/dy=\frac{\Gamma(\theta+1)}{\Gamma(1/2)\Gamma(\theta
+1/2)}\frac
{pqy^{\theta-1/2}{(1-y)}^{\theta
-1/2}}{{(p^{2}(1-y)+{q}^{2}y)}^{1+\theta}}
\]
for $0<y<1$. We also obtain results for time spent positive on
certain random subsets of $[0,1]$, and also develop some
interesting stochastic equations. As a highlight, we obtain
explicit results for the case of $\theta=\alpha$, corresponding to
the the time spent positive of a Bessel bridge on $[0,1]$. In this
case the best previous expressions were obtained independently in
\cite{JLP,Yano} (see also \cite{LijoiMean}). For some other related
works see \cite{BerYor2,Kas,KasYano,KotaniWatanabe,Watanabe,Yanos}.

Consider now the following stochastic equations and Cauchy--Stieltjes
transforms that can be found in \cite{JLP} with further
references; for $\theta>0$,
%
%
\begin{equation}\label{PDstoch1}
P_{\alpha,\theta}(p)=\beta_{\theta,1}P_{\alpha,\theta
}(p)+(1-\beta
_{\theta,1})P_{\alpha,0}(p)
\end{equation}
and for $\theta>-\alpha$,
%
%
\begin{equation}\label{PDstoch2}
P_{\alpha,\theta}(p)=\beta_{\theta+\alpha,1-\alpha}P_{\alpha
,\theta
+\alpha}(p)+(1-\beta_{\theta+\alpha,1-\alpha})\xi_{p}.
\end{equation}
Additionally there are the Cauchy--Stieltjes transforms for
$\theta>0$,
%
%
\begin{equation} \label{PCS1}
\mathcal{C}_{\theta}(\lambda;P_{\alpha,\theta
}(p))={{\bigl(q+(1+{\lambda
)}^{\alpha}p\bigr)}^{-{\theta}/{\alpha}}}
=e^{-\theta\psi^{(p)}_{\alpha,0}(\lambda)},
\end{equation}
where
\[
\psi^{(p)}_{\alpha,0}(\lambda)=\mathbb{E}\bigl[\log\bigl(1+\lambda
P_{\alpha,0}(p)\bigr)\bigr]=\mathbb{E}^{(p)}_{\alpha,0}[\log(1+\lambda
A^{+}_{1})]
\]
and for $\theta>-\alpha$,
%
%
\begin{eqnarray}\label{PCS2}
\mathcal{C}_{1+\theta}(\lambda;P_{\alpha,\theta}(p))&=&\frac
{{(1+\lambda
)}^{\alpha-1}p+(1-p)}
{{(q+(1+{\lambda)}^{\alpha}p)}^{({\theta+\alpha})/{\alpha}}}
\nonumber\\[-8pt]\\[-8pt]
&=&
\mathcal{C}_{\theta+\alpha}(\lambda;P_{\alpha,\theta+\alpha
}(p))\mathcal
{C}_{1-\alpha}(\lambda;\xi_{p}).\nonumber
\end{eqnarray}

The first equation (\ref{PDstoch1}) shows that
$P_{\alpha,\theta}(p)\stackrel{d}{=}M_{\theta}(F_{P_{\alpha,0}(p)})$
for $\theta>0$. The second equation (\ref{PDstoch2}) (see, for
instance, \cite{BertoinGoldschmidt2004,Dong2006,IJ2001,Pit06}
for some other interpretations and
applications) can be traced to Pitman and Yor (\cite{PY92},
Theorem 1.3.1) and Perman, Pitman and Yor (\cite{PPY92}, Theorem
3.8, Lemma 3.11) as follows; Let
$A^{+}_{G_{1}}=\int_{0}^{G_{1}}\indic_{(B_{s}>0)}\,ds$ denote the
time spent\vspace*{-1pt} positive of $\mathcal{B}$ up till time $G_{1}$ which
is the time of the last zero of $\mathcal{B}$ before time $1$.
Then under $\mathbb{P}^{(p)}_{\alpha,\theta}$, there is the
equivalence
\[
(A^{+}_{G_{1}},G_{1})\stackrel{d}{=}\bigl(G_{1}A^{(br)}_{1},G_{1}\bigr)\stackrel
{d}{=}(\beta_{\theta+\alpha,1-\alpha}P_{\alpha,\theta+\alpha
}(p),\beta
_{\theta+\alpha,1-\alpha}).
\]
This shows that (\ref{PDstoch2}) can be rewritten in terms of the
following decomposition:
\[
A^{+}_{1}\stackrel{d}{=}A^{+}_{G_{1}}+(1-G_{1})\xi_{p}\stackrel
{d}{=}G_{1}A^{(br)}_{1}+(1-G_{1})\xi_{p}.
\]
See, for example, Enriquez, Lucas and Simenhaus \cite{Enriquez} for
an interesting recent application of this expression.

We now show that the density of $P_{\alpha,\theta}(p)$ can be
expressed in terms of the density of $X_{\alpha,\theta}$.
Hereafter, define
\[
r_{p}(y)=\frac{y}{(c(1-y))}
\]
for
$c^{\alpha}=p/(1-p)$.
\begin{prop}\label{prop51}
For $\theta>-\alpha$, let
$R_{\alpha,\theta}=cX_{\alpha,\theta}/(cX_{\alpha,\theta}+1)$.
Then
\[
\mathbb{P}^{(p)}_{\alpha,\theta}(A^{+}_{1}\in
dy)=\frac{{(1-y)}^{\theta}}{{(1-p)}^{{\theta}/{\alpha
}}}\mathbb
{P}(R_{\alpha,\theta}\in
dy),
\]
where $A^{+}_{1}\stackrel{d}{=}P_{\alpha,\theta}(p)$. Hence as
special cases, using Propositions \ref{prop25} and \ref{prop26}:
\begin{longlist}
\item $\mathbb{P}^{(p)}_{\alpha,1}(A^{+}_{1}\in
dy)/dy={(1-y)}^{-1}p^{-1/\alpha}\Delta_{\alpha,1}(r_{p}(y))=\Omega
_{\alpha,1}(y)$.
\item For $0<\theta\leq1-\alpha$, and
$\sigma^{*}=\theta+\alpha$,
\[
\mathbb{P}^{(p)}_{\alpha,\theta}(A^{+}_{1}\in
dy)/dy=\frac{\theta{(1-y)}^{\theta-2}}{{(1-p)}^{({\theta
-1})/{\alpha
}}p^{{1}/{\alpha}}}
\int_{0}^{1}
\frac{\Delta_{\alpha,\sigma^{*}}(r_{p}(y)/u)}{u{(1-u)}^{1-\theta}}\,du,
\]
where $\Delta_{\alpha,\sigma^{*}}(x)\ge0$ is the density of
$X^{(\sigma^{*})}_{\alpha,1}$.
\end{longlist}
\end{prop}
\begin{pf}
From (\ref{subX}) it follows that, for measurable
functions $g$,
\[
\mathbb{E}[g(cX_{\alpha,\theta})]=\frac{{(1-p)}^{\theta/\alpha
}}{\mathbb
{E}[S^{-\theta}_{\alpha}]}\mathbb{E}^{(p)}_{\alpha,0}\biggl[g\biggl(
\frac{A^{+}_{\tau_{1}}}{A^{-}_{\tau_{1}}}\biggr)(A^{-}_{\tau
_{1}})^{-\theta}\biggr].
\]
The result is concluded by showing that
\[
\mathbb{E}^{(p)}_{\alpha,\theta}[g(A^{+}_{1})]:=\frac{1}{\mathbb
{E}[S^{-\theta}_{\alpha}]}
\mathbb{E}^{(p)}_{\alpha,0}[g(A^{+}_{\tau_{1}}/\tau_{1})\tau
^{-\theta}_{1}]
\]
is equal to ${{(1-p)}^{-\theta/\alpha}}
{\mathbb{E}[g(R_{\alpha,\theta}){(1+X_{\alpha,\theta})}^{-\theta}]}.
$ But this follows from
$\tau_{1}=(A^{+}_{\tau_{1}}/A^{-}_{\tau_{1}}+1)A^{-}_{\tau_{1}}$.
\end{pf}

Pitman and Yor (\cite{PY97length}, Proposition 15) establish an
interesting relationship between the densities of $A^{+}_{1}$ and
$A^{+}_{G_{1}}$ under the law $\mathbb{P}^{(p)}_{\alpha,0}$.
Making no changes to the essence of their clever argument, one can
easily extend this result to all $(\alpha,\theta)$. Combining this
with Proposition \ref{prop51} yields the relationships, for $0<p<1$,
%
%
\begin{eqnarray}\label{Denrelation}
\mathbb{P}^{(p)}_{\alpha,\theta}(A^{+}_{G_{1}}\in
dy) &=&
\frac{1-y}{(1-p)}\mathbb{P}^{(p)}_{\alpha,\theta}(A^{+}_{1}\in dy)
\nonumber\\[-8pt]\\[-8pt]
&=& \frac{{(1-y)}^{1+\theta}}{{(1-p)}^{({\theta+\alpha})/{\alpha
}}}\mathbb{P}(R_{\alpha,\theta}\in
dy).\nonumber
\end{eqnarray}
Recall that under $P^{(p)}_{\alpha,\theta}$,
$A^{+}_{G_{1}}\stackrel{d}{=}\beta_{\theta+\alpha,1-\alpha
}P_{\alpha
,\theta+\alpha}(p)$.
The next result describes interesting properties of
generalizations of this variable.
\begin{prop}\label{prop52}
For $\tau>0$ and $0<\sigma\leq1$, let
$V\stackrel{d}{=}\beta_{({\tau\sigma}/{\alpha},{\tau
(1-\sigma
)}/{\alpha})}$
and hence is a Dirichlet mean satisfying the stochastic equation
\[
V\stackrel{d}{=}\beta_{\tau/\alpha,1}V'+(1-\beta_{\tau/\alpha
,1})\xi
_{\sigma}
\]
for $V\stackrel{d}{=}V'$. Then for $0<p\leq1$,
%
%
\begin{equation}\label{betaP}
\beta_{\tau\sigma,\tau(1-\sigma)}P_{\alpha,\tau\sigma
}(p)\stackrel
{d}{=}P_{\alpha,\tau}(pV)=M_{\tau}\bigl(F_{P_{\alpha,0}(p\xi_{\sigma})}\bigr).
\end{equation}
This leads to the stochastic equations, for $0<p\leq1$,
%
%
\begin{eqnarray}\label{PV1}
P_{\alpha,\tau}(pV) &\stackrel{d}{=}& \beta_{\tau,1}P'_{\alpha
,\tau
}(pV')+(1-\beta_{\tau,1})P_{\alpha,0}(p\xi_{\sigma})
\nonumber\\[-8pt]\\[-8pt]
&\stackrel{d}{=}&
\beta_{\tau,1}P'_{\alpha,\tau}(pV)+(1-\beta_{\tau,1})P_{\alpha,0}(pV).\nonumber
\end{eqnarray}
In the first expression $P'_{\alpha,\tau}(pV')$ denotes a random
variable equivalent only in distribution to $P_{\alpha,\tau}(pV)$.
However in the second equation $V$ is the same variable.
\end{prop}
\begin{pf}
First note that\vspace*{-2pt}
$\beta_{\tau\sigma,\tau(1-\sigma)}P_{\alpha,\tau\sigma
}(p)\stackrel
{d}{=}M_{\tau}(F_{P_{\alpha,0}(p\xi_{\sigma})})$,
follows from (\ref{betascaling}). Note also by the definition of
$P_{\alpha,0}(p)$ it is easy to see that
$P_{\alpha,0}(p)\xi_{\sigma}\stackrel{d}{=}P_{\alpha,0}(p\xi
_{\sigma})$.
In order to establish the rest of (\ref{betaP}) we can check
Cauchy--Stieltjes transforms of order $\tau$ using (\ref{PCS1}).
For the variable appearing on the left of (\ref{betaP}) this is
easy to calculate. Applying this to $P_{\alpha,\tau}(pV)$
conditioned on $V$, its final evaluation rests on the simple
equality
\[
\bigl(1-pV+(1+\lambda)^{\alpha}pV\bigr)=1+[(1+\lambda)^{\alpha}-1]pV.
\]
Taking the Cauchy--Stieltjes transform of order $\tau/\alpha$ for
$V$ yields the result. The second equality in (\ref{PV1}) is then
due to (\ref{PDstoch1}).
\end{pf}

Next is one of our main distributional results which is an
analogue of Proposition \ref{prop26} but also highlights the role of
various randomly skewed processes.
\begin{prop}\label{prop53}
For $0<\sigma\leq1$, set
$R^{(\sigma)}_{\alpha,1}=cX^{(\sigma)}_{\alpha,1}/(cX^{(\sigma
)}_{\alpha,1}+1)$,
then the density of
\[
\beta_{\sigma,(1-\sigma)}P_{\alpha,\sigma}(p)\stackrel
{d}{=}P_{\alpha
,1}\bigl(p\beta_{({\sigma}/{\alpha},{(1-\sigma)}/{\alpha})}\bigr)
\]
is for $0<y<1$, equivalent to
${(1-p)}^{-{\sigma}/{\alpha}}(1-y)\mathbb{P}(R^{(\sigma
)}_{\alpha
,1}\in
dy)/dy$, and is given explicitly as
%
%
\begin{equation}\label{denB}\qquad
\Omega_{\alpha,\sigma}(y)=\frac{1}{\pi}
\frac{y^{\sigma-1}\sin(\rho_{\alpha,\sigma
}({[r_{p}(y)]}^{\alpha
}) )}
{{[y^{2\alpha}q^{2}+2qpy^{\alpha}{(1-y)}^{\alpha}\cos(\alpha
\pi)+{(1-y)}^{2\alpha}p^{2}]}^{{\sigma}/({2\alpha})}},
\end{equation}
where $\rho_{\alpha,\sigma}$ is defined in (\ref{short}). In
particular $\Omega_{\alpha,1}(y)$ is the density of
$P_{\alpha,1}(p)$ and $\Omega_{\alpha,\alpha}(y)$ is the density
of $\beta_{\alpha,1-\alpha}P_{\alpha,\alpha}(p)$. In addition:
\begin{longlist}
\item
if $0\leq\theta\leq1-\alpha$, then for
$\sigma^{*}=\theta+\alpha$,
\[
\beta_{\theta+\alpha,1-\alpha}P_{\alpha,\theta+\alpha
}(p)\stackrel
{d}{=}P_{\alpha,1+\theta}\bigl(p\beta_{({\sigma^{*}}/{\alpha
},({1-\alpha})/{\alpha})}\bigr)\stackrel{d}{=}\beta_{1,\theta}P_{\alpha
,1}\bigl(p\beta
_{({\sigma^{*}}/{\alpha},({1-\sigma^{*}})/{\alpha})}\bigr)
\]
and, using (\ref{Denrelation}), there is the explicit formula
determining the densities of $A^{+}_{1}$ and $A^{+}_{G_{1}}$,
%
%
\begin{eqnarray}
\label{coagden}
\mathbb{P}^{(p)}_{\alpha,\theta}(A^{+}_{G_{1}}\in
dy)/dy &=&
\frac{1-y}{(1-p)}\mathbb{P}^{(p)}_{\alpha,\theta}(A^{+}_{1}\in
dy)/dy
\nonumber\\[-8pt]\\[-8pt]
&=&
\int_{0}^{1}\frac{\theta\Omega_{\alpha,\sigma
^{*}}(y/u)}{u{(1-u)}^{1-\theta}}\,du.\nonumber
\end{eqnarray}
When $\theta=0$, (\ref{coagden}) is $\Omega_{\alpha,\alpha}(y)$
which agrees with Pitman and Yor (\cite{PY97length}, Proposition~15).
\item As in Proposition \ref{prop27}, for any $\theta>0$, set
$\theta
=\sum_{j=1}^{k}\theta_{j}$ for some integer $k$ and $\theta_{j}>0$.
This leads to the representation
\[
P_{\alpha,\theta}(p)\stackrel{d}{=}\sum_{j=1}^{k}D_{j}P_{\alpha
,\theta_{j}}(p)
\]
for independent variables $P_{\alpha,\theta_{j}}(p)$ and
$(D_{1},\ldots,D_{k}$) a Dirichlet vector as in
Proposition \ref{prop27}. When $\theta_{j}$ are chosen such that $0<
\theta_{j}\leq1-\alpha$ each variable has
$\mathbb{P}^{(p)}_{\alpha,\theta_{j}}(A^{+}_{1}\in dy)$ given by
(\ref{coagden}) with $\sigma^{*}=\theta_{j}+\alpha$. It suffices
to choose $\theta_{j}=\theta/k$, for $0<\theta\leq k(1-\alpha)$.
If $\theta=k$, one may set $\theta_{j}=1$ and use
$\Omega_{\alpha,1}$.
\end{longlist}
\end{prop}
\begin{pf}
The various representations of the random variables
are due to Proposition \ref{prop52} and otherwise an application of the
beta/gamma calculus. The density $\Omega_{\alpha,\sigma}$ is
obtained from $\Delta_{\alpha,\sigma}$, which is justified by the
exponential tilting relationships discussed in
James (\cite{JamesBernoulli}, Section \ref{sec3}; see
also \cite{JamesGamma}).
\end{pf}
\begin{rem}\label{rskewed} The random variable $P_{\alpha,\tau}(pV)$
described
in Proposition \ref{prop52} has law,
\[
\mathbb{P}\bigl(P_{\alpha,\tau}(pV)\in
dx\bigr)=\mathbb{P}^{(pV)}_{\alpha,\tau}(A^{+}_{1}\in
dx):=\int_{0}^{1}\mathbb{P}^{(pu)}_{\alpha,\tau}(A^{+}_{1}\in
dx)f_{V}(u)\,du.
\]
That is, it may be read as the time spent positive up till one of
a process $\mathcal{B}$ whose excursion lengths, conditional on
$V$, follow a $\operatorname{PD}(\alpha,\tau)$ distribution and is otherwise
randomly skewed by $pV$. See also Aldous and
Pitman (\cite{AldousT}, Section \ref{sec51}) for connections with
$T$-partitions. This is made clear, as follows; For $(L_{t}; 0\leq
t\leq1)$ governed by $\operatorname{PD}(\alpha,\theta)$, and letting
$\bar{L}_{t}=L_{t}/L_{1}$, there is the equivalence
\[
P_{\alpha,\theta}(u)\stackrel{d}{=}\inf\{t\dvtx\bar{L}_{t}\ge
u\},\qquad 0\leq u\leq1.
\]
In other\vspace*{-2pt} words, letting $P^{(-1)}_{\alpha,\theta}(\cdot)$ denote
the random \textit{quantile} function of $P_{\alpha,\theta}$, it
follows that $\bar{L}_{t}\stackrel{d}{=}P^{(-1)}_{\alpha,\theta}(t),
0\leq t\leq1$. See the next section, Section \ref{sec6}, for more general
$V$.
\end{rem}
\begin{rem}\label{remark52}
In reference to Propositions \ref{prop52} and \ref
{prop53}, setting $Q_{\alpha,\tau}(\sigma,p)=\beta_{\tau\sigma,\tau
(1-\sigma)}P_{\alpha,\tau\sigma}(p)$
leads to a well-defined bivariate process
$(Q_{\alpha,\tau}(\sigma,p)\dvtx 0\leq p\leq1, 0<\sigma\leq1)$, that
has some natural connections to the coagulation operations
discussed in Pitman \cite{Pit99}. This observation may be deduced
from the subordinator representation given in Pitman and
Yor (\cite{PY97}, Proposition 21). When $p=1$,
$Q_{\alpha,\tau}(\sigma,1)$ is a Dirichlet process which
corresponds to the operation of coagulating
$\operatorname{PD}(\alpha,\tau)$ by $\operatorname{PD}(0,\tau/\alpha)$. In
general one may write
\[
Q_{\alpha,\tau}(\sigma,p)\stackrel{d}{=}P_{0,\tau}(\sigma
)P_{\alpha,\tau
\sigma}(p)\stackrel{d}{=}P_{\alpha,\tau}(pP_{0,\frac{\tau}{\alpha
}}(\sigma)).
\]
We will not elaborate on this here except to note that in
connection with results for the standard $U$-coalescent, setting
$\tau=\alpha$, $p=1$ one recovers (\cite{Pit99}, Corollary 33, and
Proposition 32). In a distributional sense, that is, without the
nice interpretation, one also recovers (\cite{Pit99}, Corollary
16) by setting $\sigma=\alpha=e^{-t}$ and $\tau=1,p=1$. When
$p\neq1$, we expect that one can obtain new, but related,
interpretations of $Q_{\alpha,\tau}(\sigma,p)$.
\end{rem}

\subsection{Some special cases}\label{sec51}

Note that, as in \cite{JLP} (combined with Proposition~\ref{prop52}), one
can rewrite (\ref{PDstoch2}) as
\begin{eqnarray*}
P_{\alpha,\theta}(p)&\stackrel{d}{=}&\beta_{\theta+\alpha
,1-\alpha
}P_{\alpha,\theta+\alpha}(p)(1-\xi_{p})+\bigl(1-\beta_{\theta+\alpha
,1-\alpha
}P_{\alpha,\theta+\alpha}(q)\bigr)
\xi_{p}\\
&\stackrel{d}{=}&
P_{\alpha,1+\theta}\bigl(p\beta_{(({\theta+\alpha})/{\alpha},
({1-\alpha})/{\alpha})}\bigr)(1-\xi_{p})\\
&&{} + \bigl(1-
P_{\alpha,1+\theta}\bigl(q\beta_{(({\theta+\alpha})/{\alpha},
({1-\alpha})/{\alpha})}\bigr)\bigr)
\xi_{p}.
\end{eqnarray*}
Besides giving an alternate mixture representation in terms of
easily interpreted random variables, this also suggests that one
can obtain the density of $P_{\alpha,\theta}(p)$ if one knows the
density of $P_{\alpha,\theta+\alpha}(p)$. In \cite{JLP}, it was
noted that this could be applied for $\theta+\alpha=1$, which
yields an expression for the density of $P_{\alpha,1-\alpha}(p)$.
In view of Proposition \ref{prop53} we see that such a density
representation can be extended to any $0\leq\theta\leq1-\alpha$.
Of course in terms of a density representation this is not as good
as the expression one can obtain from (\ref{coagden}), since it
would have to be used twice. In this section we look at some
specific cases of random variables that have either appeared in
the literature or we anticipate might be of some interest.
\begin{example}[{[$(\alpha,1-\alpha)$ a distribution relevant to
phylogenetic models]}]\label{example51}
As noted in \cite{JLP} the case of $P_{\alpha,1-\alpha}(p)$
equates in distribution to the limit of a phylogenetic tree model
appearing in (\cite{Haas}, Proposition 20). Here using
Proposition~\ref{prop53} we obtain a slight improvement over the
density given in (\cite{JLP}, Corollary 6.1). Since under
$\mathbb{P}^{(p)}_{\alpha,1-\alpha}$,
$A^{+}_{G_{1}}\stackrel{d}{=}\beta_{1,1-\alpha}P_{\alpha,1}(p)$, we
have
%
%
\begin{eqnarray}
\label{phyloden}
\mathbb{P}^{(p)}_{\alpha,1-\alpha}(A^{+}_{G_{1}}\in
dy)/dy &=&
\frac{1-y}{(1-p)}\mathbb{P}^{(p)}_{\alpha,1-\alpha}(A^{+}_{1}\in
dy)/dy
\nonumber\\[-8pt]\\[-8pt]
&=&
\int_{0}^{1}\frac{(1-\alpha)\Omega_{\alpha
,1}(y/u)}{u{(1-u)}^{\alpha}}\,du.\nonumber
\end{eqnarray}
\end{example}
\begin{example}[{[The case of
$\beta_{1+\alpha,1-\alpha}P_{\alpha,1+\alpha}(p)$]}]\label{example52}
Under $\mathbb{P}^{(p)}_{\alpha,1}$,
\[
A^{+}_{G_{1}}\stackrel{d}{=}\beta_{1+\alpha,1-\alpha}P_{\alpha
,1+\alpha
}(p)\stackrel{d}{=}P_{\alpha,2}\bigl(p\beta_{(({1+\alpha})/{\alpha
},({1-\alpha})/{\alpha})}\bigr).
\]
Hence its density is given by
\[
\mathbb{P}^{(p)}_{\alpha,1}(A^{+}_{G_{1}}\in
dy)/dy=\frac{(1-y)}{(1-p)}\Omega_{\alpha,1}(y).
\]
In view of the literature related to Section \ref{sec4} we believe this
variable will be of interest.
\end{example}

We now address some harder cases.
\begin{example}[{[$(\alpha,\alpha)$ occupation time of a Bessel bridge]}]\label{example53}
Obtaining density expressions for the general case of $A^{+}_{1}$
when $\mathcal{B}$ is a Bessel bridge has been difficult, except
for the case of $\alpha=1/2$. Due to the importance of the
$\operatorname{PD}(\alpha,\alpha)$ family this quantity arises in many contexts (see,
for instance, Aldous and Pitman \cite{AldousT}). The best
results were obtained independently by Yano \cite{Yano} and James,
Lijoi and Pr\"{u}nster
\cite{JLP}, who give expressions in terms of Abel-type transforms, that
is to say, integrals of possibly nonnegative functions. Hence
this does not yield a mixture representation for
$P_{\alpha,\alpha}(p)\stackrel{d}{=}A^{+}_{1}$. Here we show how our
results in the previous section can be used to achieve this. Under
$\mathbb{P}^{(p)}_{\alpha,\alpha}$,
\[
A^{+}_{G_{1}}\stackrel{d}{=}\beta_{2\alpha,1-\alpha}P_{\alpha
,2\alpha
}(p)\stackrel{d}{=}P_{\alpha,1+\alpha}\bigl(p\beta_{(2,({1-\alpha
})/{\alpha})}\bigr).
\]
Hence for $\alpha\leq1/2$, we can apply statement (i) of
Proposition \ref{prop53} writing
\[
A^{+}_{G_{1}}\stackrel{d}{=}\beta_{1,\alpha}P_{\alpha,1}\bigl(p\beta
_{(2,({1-2\alpha})/{\alpha})}\bigr)
\]
to get
%
%
\begin{eqnarray}
\label{besbridge}
\mathbb{P}^{(p)}_{\alpha,\alpha}(A^{+}_{G_{1}}\in
dy)/dy &=&
\frac{1-y}{(1-p)}\mathbb{P}^{(p)}_{\alpha,\alpha}(A^{+}_{1}\in
dy)/dy
\nonumber\\[-8pt]\\[-8pt]
&=&
\int_{0}^{1}\frac{\alpha\Omega_{\alpha,2\alpha
}(y/u)}{u{(1-u)}^{1-\alpha}}\,du,\nonumber
\end{eqnarray}
where
\[
\Omega_{\alpha,2\alpha}(y)=\frac{2\sin(\pi
\alpha)}{\pi}\frac{py^{2\alpha-1}(1-y)^{\alpha}[qy^{\alpha}+\cos
(\pi
\alpha)p(1-y)^{\alpha}]
}{{[y^{2\alpha}q^{2}+2qpy^{\alpha}{(1-y)}^{\alpha}\cos(\pi\alpha
)+{(1-y)}^{2\alpha}p^{2}]}^{2}}
\]
is the density of
$P_{\alpha,1}(p\beta_{(2,({1-2\alpha})/{\alpha})})$.

When $\alpha>1/2$, we, at present, need to resort to statement
(ii) of Proposition~\ref{prop53}. So, for instance, for $\alpha\leq2/3$,
it follows that
\[
P_{\alpha,\alpha}(p)\stackrel{d}{=}\beta_{\alpha/2,\alpha
/2}P_{\alpha
,\alpha/2}(p)+(1-\beta_{\alpha/2,\alpha/2})P'_{\alpha,\alpha/2}(p),
\]
where $P_{\alpha,\alpha/2}(p)$ and $P'_{\alpha,\alpha/2}(p)$ are
i.i.d. variables having distribution\break
$\mathbb{P}^{(p)}_{\alpha,\alpha/2}(A^{+}_{1}\in dy)$ obtainable
from (\ref{coagden}).
\end{example}
\begin{example}[{[$(\alpha,\alpha-1)$, and fragmentation equations]}]\label{example54}
Suppose that we are interested in the case of $\theta=\alpha-1$,
that under $\mathbb{P}^{(p)}_{\alpha,\alpha-1}$,
$P_{\alpha,\alpha-1}(p)\stackrel{d}{=}A^{+}_{1}$. Of course this only
makes sense for $\alpha>1/2$. Notice that
\[
A^{+}_{G_{1}}\stackrel{d}{=}\beta_{2\alpha-1,1-\alpha}P_{\alpha
,2\alpha-1}(p),
\]
so we can apply (\ref{coagden}) directly if $\alpha\leq2/3$. It
is then interesting to note what other quantities we can obtain.
We can use another stochastic equation that takes the form
\[
P_{\alpha,\theta}(p)\stackrel{d}{=}\beta_{\theta+\alpha\delta
,1-\alpha
\delta}P_{\alpha,\theta+\alpha\delta}(p)+
(1-\beta_{\theta+\alpha\delta,1-\alpha\delta})P_{\alpha,-\alpha
\delta}(p)
\]
for $\theta>-\alpha\delta$, $0<\delta\leq1$. Note that this
equation is not well known but it is simple to check. Furthermore,
a close inspection shows that it is a nice way to code
Pitman's \cite{Pit99} fragmentation. As as special case, set
$\delta=(1-\alpha)/\alpha$ and $\theta=\alpha$ to obtain
\[
P_{\alpha,\alpha}(p)\stackrel{d}{=}\beta_{1,\alpha}P_{\alpha
,1}(p)+(1-\beta_{1,\alpha})P_{\alpha,\alpha-1}(p).
\]
\end{example}

\section{Power scaling property and randomly skewed processes}\label{sec6}

We saw that in the previous section the random quantity
$P_{\alpha,\tau}(pV)$, where $V$ is a beta variable, occurs
naturally and plays an interesting role. An interpretation in
terms of occupation times of randomly skewed processes is
mentioned in Remark \ref{rskewed}, and an interpretation via coagulation
processes is hinted at in Remark \ref{remark52}. Also there is the surprising
stochastic equation in (\ref{PV1}). There is also a related result
given in Proposition \ref{prop22}. One may wonder if properties of this
sort only hold for beta random variables. We show in the next
result, which was first obtained in \cite{JamesGamma}, that there
is a considerable generalization.
\begin{prop}\label{prop61} Let $\mathcal{R}\stackrel{d}{=}M_{\tau/\alpha}(F_{R})$ and
$\mathcal{Q}\stackrel{d}{=}M_{\tau/\alpha}(F_{Q})$ denote Dirichlet
means with parameters $(\tau/\alpha, R)$ and $(\tau/\alpha, Q)$
where $R$ is a nonnegative random variable, and $Q$ is a random
variable taking values in $[0,1]$. Equivalently,
%
%
\begin{eqnarray}\label{Qmeans}
\mathcal{R}&\stackrel{d}{=}&\beta_{({\tau}/{\alpha},1)}\mathcal
{R}+\bigl(1-\beta_{({\tau}/{\alpha},1)}\bigr)R
\quad\mbox{and}\nonumber\\[-8pt]\\[-8pt]
\mathcal{Q}&\stackrel{d}{=}&\beta_{({\tau}/{\alpha},1)}\mathcal
{Q}+\bigl(1-\beta_{({\tau}/{\alpha},1)}\bigr)Q,\nonumber
\end{eqnarray}
which implies that $\mathcal{Q}\stackrel{d}{=}M_{\tau/\alpha}(F_{Q})$
takes it values in $[0,1]$. If $Q$ is a constant, then
$M_{\tau/\alpha}(F_{Q})=Q$. Then the following results hold:
\begin{longlist}
\item $\mathcal{R}^{1/\alpha}X_{\alpha,\tau}\stackrel
{d}{=}M_{\tau
}(F_{X_{\alpha}R^{1/\alpha}})$,
that is,
\[
\mathcal{R}^{1/\alpha}X_{\alpha,\tau}\stackrel{d}{=}\beta_{\tau
,1}\mathcal{R}^{1/\alpha}X_{\alpha,\tau}+(1-\beta_{\tau
,1})R^{1/\alpha
}X_{\alpha}.
\]
\item $P_{\alpha,\tau}(\mathcal{Q})$ is a Dirichlet mean with
parameters $(\tau,P_{\alpha,0}(Q))$, and satisfies,
%
%
\begin{eqnarray}\label{PQ}
P_{\alpha,\tau}(\mathcal{Q}) &\stackrel{d}{=}& \beta_{\tau
,1}P'_{\alpha
,\tau}(\mathcal{Q}')+(1-\beta_{\tau,1})P_{\alpha,0}(Q)
\nonumber\\[-8pt]\\[-8pt]
&\stackrel{d}{=}&
\beta_{\tau,1}P'_{\alpha,\tau}(\mathcal{Q})+(1-\beta_{\tau
,1})P_{\alpha
,0}(\mathcal{Q}),\nonumber
\end{eqnarray}
where $P'_{\alpha,\tau}(\mathcal{Q}')$ is equivalent only in
distribution to $P_{\alpha,\tau}(\mathcal{Q})$, but $\mathcal{Q}$
is the
same variable.
\end{longlist}
\end{prop}
\begin{pf}
The first result follows by noting
\[
\tau
\mathbb{E}[\log(1+\lambda
X_{\alpha}R^{1/\alpha})]=(\tau/\alpha)\mathbb{E}[\log(1+\lambda
^{\alpha}R)],
\]
which gives the negative logarithm of $\mathcal{C}_{\tau}(\lambda;
\mathcal{R}^{1/\alpha}X_{\alpha,\tau})=\mathcal{C}_{\tau/\alpha
}(\lambda
^{\alpha};
\mathcal{R})$. For the second, evaluate
$\mathcal{C}_{\tau}(\lambda, P_{\alpha,\tau}(\mathcal{Q}))$
conditional on $\mathcal{Q}$ and then notice similar to
Proposition \ref{prop52}, that the transform of order $\tau$ coincides with
the exponential of
\[
\log\mathcal{C}_{\tau/\alpha}\bigl({(1+\lambda)}^{\alpha}-1;
\mathcal{Q}\bigr)=-\frac{\tau}{\alpha}\mathbb{E}\bigl[\psi^{(Q)}_{\alpha,0}\bigr].
\]
It remains then to apply (\ref{PDstoch1}) to get the second
equality in (\ref{PQ}).
\end{pf}
\begin{rem}
Notice that setting $R=\xi_{\sigma}$ and $Q=p\xi_{\sigma}$ we
recover Propositions \ref{prop22} and \ref{prop52}. Setting
$R=X_{\delta}$ for $0<\delta<1$ leads to the identity
$X^{1/\alpha}_{\delta,\tau/\alpha}X_{\alpha,\tau}\stackrel
{d}{=}X_{\alpha\delta,\tau}$, 
since it follows from known properties of stable random variables
that
$X^{1/\alpha}_{\delta}X_{\alpha}\stackrel{d}{=}X_{\alpha\delta}$.
Furthermore, if one chooses $\mathcal{Q}:=\mathcal{Q}(u)$ such for
each fixed $u$ it satisfies (\ref{Qmeans}), and for $0<u<1$ it is
an exchangeable bridge, that is a random cumulative distribution
function, then $P_{\alpha,\tau}(\mathcal{Q}(u))$ identifies a
coagulation operation as described in Pitman (\cite{Pit06}, Lemma
5.18) (see also Bertoin \cite{BerFrag}). In particular, one
recovers Pitman's \cite{Pit99} coagulation as follows. Setting
$\mathcal{Q}=P_{\beta,\tau/\alpha}(u)$, means that
$Q=P_{\delta,0}(u)$, leading easily to,
$P_{\alpha,0}(P_{\delta,0}(u))\stackrel{d}{=}P_{\alpha\delta,0}(u)$,
which implies
\[
P_{\alpha,\tau}(P_{\delta,\tau/\alpha}(u))\stackrel
{d}{=}P_{\alpha\delta
,\tau}(u).
\]
We shall discuss other applications of Proposition \ref{prop61} and related
identities elsewhere.
\end{rem}

\section{Subordinators and symmetric generalized Linnik laws and processes}\label{sec7}

Using Proposition \ref{prop61}, we define processes $(T_{\alpha}(\tau),
\tau\ge0)$ and $(\hat{T}_{\alpha}(\tau), \tau\ge0)$, such that
for each fixed $\tau>0$,
\[
T_{\alpha}(\tau)\stackrel{d}{=}\gamma_{\tau}\mathcal{R}^{1/\alpha
}X_{\alpha,\tau}
\stackrel{d}{=}\chi_{\alpha,\tau}\mathcal{R}^{1/\alpha}
\stackrel{d}{=}\gamma_{\tau}M_{\tau}(F_{X_{\alpha}R^{1/\alpha}})
\]
and
\[
\hat{T}_{\alpha}(\tau)\stackrel{d}{=}\gamma_{\tau}P_{\alpha,\tau
}(\mathcal{Q})
\]
are $\operatorname{GGC}(\tau, X_{\alpha}R^{1/\alpha})$ and
$\operatorname{GGC}(\tau,P_{\alpha,0}(Q))$ variables, respectively. Where we are
suppressing the fact that both $\mathcal{R}$ and $\mathcal{Q}$
depend on $(\alpha,\tau)$. In fact $T_{\alpha}(\tau)$ and
$\hat{T}_{\alpha}(\tau)$ are GGC subordinators varying in
$\tau>0$. Let $\mathcal{S}_{\alpha}(t)$ denote a positive stable
subordinator such that
$\mathcal{S}_{\alpha}(1)\stackrel{d}{=}S_{\alpha}$, and let
$\hat{\mathcal{S}}_{\alpha}(t)$ denote the subordinator with
\[
-\log\mathbb{E}\bigl[e^{-\lambda
\hat{\mathcal{S}}_{\alpha}(t)}\bigr]=t[(1+\lambda)^{\alpha}-1],
\]
so that
${\hat{\mathcal{S}}_{\alpha}(1)}\stackrel{d}{=}\hat{S}_{\alpha}$
is a
random variable with density $e^{-(t-1)}f_{\alpha}(t)$. It
follows that
\[
\mathcal{S}_{\alpha}(\gamma_{{\tau}/{\alpha}}\mathcal
{R})\stackrel
{d}{=}T_{\alpha}(\tau)
\quad\mbox{and}\quad
\hat{\mathcal{S}}_{\alpha}(\gamma_{{\tau}/{\alpha}}\mathcal
{Q})\stackrel{d}{=}\hat{T}_{\alpha}(\tau).
\]
However, an important aspect of our results in Proposition \ref{prop61} is
that usage of $T_{\alpha}(\tau)$ and $\hat{T}_{\alpha}(\tau)$ does
not require working directly with the processes
$\mathcal{S}_{\alpha}$ and $\hat{\mathcal{S}}_{\alpha}$. These
L\'evy processes are attractive in terms of potential applications
arising for instance in finance, financial econometrics or
Bayesian statistics. With applications to finance in mind, it is
quite natural to use these processes as Brownian time changes
creating process $\mathcal{B}(T_{\alpha}(\cdot))$ and
$\mathcal{B}(\hat{T}_{\alpha}(\cdot))$, for $\mathcal{B}(\cdot)$
an independent Brownian motion, we take to have log characteristic
function $-\lambda^{2}$. The characteristic functions of
$\mathcal{B}(T_{\alpha}(\tau))$ and
$\mathcal{B}(\hat{T}_{\alpha}(\tau))$ can be expressed as
\[
\mathcal{C}_{{\tau}/{\alpha}}(\lambda^{2\alpha};\mathcal
{R})=e^{-{\tau}/{\alpha}\psi_{R}(\lambda^{2\alpha})}
\]
and
\[
\mathcal{C}_{{\tau}/{\alpha}}\bigl({(1+\lambda^{2})}^{\alpha
}-1;\mathcal
{Q}\bigr)=e^{-{\tau}/{\alpha}\psi_{Q}({(1+\lambda^{2})}^{\alpha}-1)}.
\]
Recall that $\mathcal{B}(\mathcal{S}_{\alpha}(\cdot))$ is a
symmetric stable process of index $(0,2]$ and
$\mathcal{B}(\hat{\mathcal{S}}_{\alpha}(\cdot))$ is a process
that includes the NIG process when $\alpha=1/2$. When $\alpha=0$
and $Q=p$, $\mathcal{B}(\hat{T}_{0}(\cdot))$, is a variance-gamma
(VG) process. The case of $\mathcal{B}(\chi_{\alpha,\theta})$,
corresponds to generalized Linnik processes considered by
Pakes \cite{Pakes} (see also \cite{DevroyeLinnik,Koz}). We
will focus on this case.

It suffices to examine the random variables
\[
\mathcal{B}(\chi_{\alpha,\theta})\stackrel{d}{=}N\sqrt{2\gamma
^{1/\alpha
}_{{\theta}/{\alpha}}S_{\alpha}}=N\sqrt{\chi_{\alpha,\theta}},
\]
where $N$ is a standard Normal random variable. For general
$\alpha$ and $\theta>0$ the extra randomization by $N$ does not
add much beyond our results for $\chi_{\alpha,\theta}$.

However, when $\alpha\leq1/2$ we are able to obtain some
interesting results which we describe below. In this case, we will
use a result of Devroye \cite{DevroyePolya} which yields a
tractable mixture representation for symmetric stable random
variables of index between $0$ and $1$.
Devroye's \cite{DevroyePolya} result is not well known but as we
shall show can be used to obtain a nice description of the density
of $\mathcal{B}_{\alpha}(T_{\alpha}(\theta))$ for all fixed
$\theta>0$.
We do, however, stress that there are many applications requiring
$\alpha>1/2$.

\subsection{$\alpha\leq12$, results based on Fejer-de la Vallee
Poussin mixtures}\label{sec71}

For symmetric stable random variables $N\sqrt
{2S_{\alpha}}$ for
$0<\alpha\leq1/2$, Devroye \cite{DevroyePolya} shows that
\[
N\sqrt{2S_{\alpha}}\stackrel{d}{=}Y/Z^{{1}/({2\alpha})},
\]
where $Y$ has a Fejer-de la Vallee Poussin density
\[
\omega(x)=\frac{1}{2\pi}{ \biggl(\frac{\sin(x/2)}{x/2}
\biggr)}^{2},\qquad
-\infty<x<\infty,
\]
and
$Z\stackrel{d}{=}\gamma_{1}(1-\xi_{2\alpha})+\gamma_{2}\xi
_{2\alpha}$.
It follows that the density of a symmetric stable of index between
$[0,1]$ is
\[
2\alpha\int_{0}^{\infty}\omega(xy)y^{2\alpha}e^{-y^{2\alpha}}
[(1-2\alpha)+2\alpha y^{2\alpha}]\,dy.
\]
Hence as a mild extension of Devroye (\cite{DevroyePolya}, Example
B), that is, the simple symmetric Linnik variable corresponding to
$\theta=\alpha$, we have for all $\theta>0$
\[
N\sqrt{2\chi_{\alpha,\theta}}\stackrel{d}{=}Y/W^{{1}/({2\alpha})}
\qquad\mbox{where }
W\stackrel{d}{=}\frac{\gamma_{1}}{\gamma_{{\theta}/{\alpha
}}}(1-\xi
_{2\alpha})+\frac{\gamma_{2}}{\gamma_{{\theta}/{\alpha}}}\xi
_{2\alpha},
\]
is a mixture of Pareto variables, having density for
$0<\alpha\leq1/2$
\[
\frac{\theta[(1+2\theta)
w+(1-2\alpha)]}{\alpha{(1+w)}^{2+\theta/\alpha}},\qquad w>0.
\]
Naturally this representation extends to all
$\mathcal{B}({T}_{\alpha}(\cdot))$, provided that $\alpha\leq
1/2$. In particular,
\[
\mathcal{B}({T}_{\alpha}(\theta))\stackrel{d}{=}Y/\tilde{W}^{
{1}/({2\alpha})}\qquad
\mbox{where now }
\tilde{W}\stackrel{d}{=}\frac{\gamma_{1}}{\gamma_{{\theta
}/{\alpha
}}\mathcal{R}}(1-\xi_{2\alpha})+\frac{\gamma_{2}}
{\gamma_{{\theta}/{\alpha}}\mathcal{R}}\xi_{2\alpha},
\]
for $\mathcal{R}\stackrel{d}{=}M_{\theta/\alpha}(F_{R})$. Quite
interestingly the density of $\tilde{W}$ only requires information
about the Laplace transform of
$\gamma_{\theta/\alpha}\mathcal{R}$. Let $\psi^{(1)}_{R}(x)$ and
$\psi^{(2)}_{R}(x)$ denote the first and second derivatives of
$\psi_{R}(x)$. Then the density of $\tilde{W}$ is given by
\[
\eta_{\alpha,\theta}(x)=\frac{\theta}{\alpha}e^{-{\theta}/{\alpha}\psi_{R}(x)}
\biggl[\psi^{(1)}_{R}(x)(1-2\alpha)+x2\alpha\biggl[\bigl(\psi^{(1)}_{R}(x)\bigr)^{2}\frac
{\theta}{\alpha}-\psi^{(2)}_{R}(x)\biggr]\biggr].
\]
From
this, we close with an interesting identity.
\begin{prop}\label{prop71} For $0<\alpha\leq1/2$, and $\theta>0$, let $V$
$\stackrel{d}{=}\beta_{1/2,1/2}$, then for $-\infty<x<\infty$
\begin{eqnarray*}
\Phi_{\alpha,\theta}(x)&=&
\mathbb{E} \biggl[\frac{|x|}{2V\chi_{\alpha,\theta}}e^{-{x^{2}}/({4V\chi_{\alpha,\theta}})} \biggr]\\ &=&
\sqrt{\frac{2}{\pi}}
\mathbb{E} \biggl[\frac{1}{|\sqrt{2\chi_{\alpha,\theta}}|}e^{-{x^{2}}/({4\chi_{\alpha,\theta}})} \biggr]\\
&=& \int_{0}^{\infty} \frac{\omega(xy)2\theta
y^{2\alpha}[(1+2\theta) y^{2\alpha}
+(1-2\alpha)]}{{(1+y^{2\alpha})}^{2+\theta/\alpha}}\,dy,
\end{eqnarray*}
which is just the density of $N\sqrt{2\chi_{\alpha,\theta}}$.
Additionally, for all fixed $\theta>0$, the density of
$\mathcal{B}({T}_{\alpha}(\theta))\stackrel{d}{=}N\sqrt{2\chi
_{\alpha
,\theta}}\mathcal{R}^{{1}/({2\alpha})}$,
satisfies
\[
\mathbb{E} \biggl[\Phi_{\alpha,\theta}\biggl(\frac{x}{\mathcal{R}^{
{1}/({2\alpha})}}\biggr)\frac{1}{\mathcal{R}^{{1}/({2\alpha})}}
\biggr]=2\alpha
\int_{0}^{\infty}\omega(xy)y^{2\alpha}\eta_{\alpha,\theta
}(y^{2\alpha})\,dy
\]
for $\mathcal{R}\stackrel{d}{=}M_{\theta/\alpha}(F_{R})$.
\end{prop}
\begin{pf}
The result follows from the derivation of the density
described above using $\omega$, in combination with a derivation
of the density based on $N^{2}\stackrel{d}{=}2\gamma_{1}V$ and
additionally Pitman and Yor (\cite{PYarcsine}, equation (29)).
\end{pf}

\section{General remark about rational case}\label{sec8}

When $\alpha$ is rational it is known (see \cite{Chaumont,Williams,ZolotarevStable})
that $S_{\alpha}$ is equivalent in distribution
to a product of independent beta and gamma variables. Extending an
argument in Chaumont and Yor (\cite{Chaumont}, pages 143 and 144) using
the gamma duplication formula, it follows that for $\alpha=m/n$
for integers, $m,n$, such that $m<n$, and all $\theta>-m/n$,
\[
{(X_{{m/n},\theta})}^{m}\stackrel{d}{=}{ \biggl(\frac{S_{
{m/n}}}{S_{{m/n},\theta}} \biggr)}^{m}
\stackrel{d}{=} \Biggl(\prod_{k=1}^{m-1}
\frac{\beta_{{\theta}/{m}+{k/n},k({1/m}-{1/n})}}
{\beta_{{k/n},k({1/m}-{1/n})}} \Biggr)
\Biggl(\prod_{k=m}^{n-1}
\frac{\gamma_{{\theta/m}+{k/n}}
}{\gamma_{{k/n}}} \Biggr),
\]
where all random variables are independent. Additionally,
\[
{\biggl(\frac{m}{S_{{m/n},\theta}} \biggr)}^{m}\stackrel{d}{=}n^{n}
\Biggl(\prod_{k=1}^{m-1}\beta_{{\theta/m}+
{k/n},k(
{1/m}-{1/n})} \Biggr)
\Biggl(\prod_{k=m}^{n-1}\gamma_{{\theta/m}+
{k/n}} \Biggr).
\]
An implication of these relationships is that one may use the
result of Springer and Thompson \cite{Springer70} to express their
densities in terms of Meijer-$G$ functions. In many cases these are
equivalent to expressions in terms of generalized Gauss
hypergeometric functions. Furthermore, it is known that Laplace
transforms of Meijer-$G$ functions are also Meijer-$G$ functions, with
known arguments. Hence, Propositions \ref{prop31} and \ref{prop51} show
that in the rational case of $\alpha=m/n$, one may express the
generalized Mittag--Leffler functions and densities for
$P_{\alpha,\theta}(p)$ in terms of Meijer-$G$ functions. Such
representations are not entirely appealing in many respects; for
instance, the density $\Delta_{m/n,\sigma}$ is a much more
desirable expression than its Meijer-$G$ counterpart. However, from
a computational viewpoint they are significant. This is due to the
fact that Meijer-$G$ functions, which constitute many special
functions, are available as built-in functions in mathematical
computational packages such as Mathematica or Maple. Naturally
many of the quantities we discussed for general $\alpha$ can be
expressed as the more general Fox-$H$ functions \cite
{MPS,Kirbook,KilbasSaigo,kir1,Kilbas04,Fox61}. However, in general,
computations for these expressions are not yet available. Hence
another contribution of our work is to give new explicit
identities for a class of Meijer-$G$ and Fox-$H$ functions. That is to
say quantities such as $\Delta_{\alpha,\sigma}$ give an explicit
form to their corresponding Fox-$H$ representation. We omit details
of this representation, but it is not difficult to obtain.


%
\printaddresses

\end{document}